\documentclass[leqno]{amsart}

\usepackage{amsmath,amssymb,amsthm}
\usepackage{mathrsfs}
\usepackage{braket}
\usepackage{graphicx}
\usepackage{amscd}
\usepackage[all]{xy}
\usepackage{comment}
\usepackage{eucal}

\theoremstyle{definition}
\newtheorem{dfn}{Definition}[section]
\newtheorem{lem}[dfn]{Lemma}
\newtheorem{cor}[dfn]{Corollary}
\newtheorem{prp}[dfn]{Proposition}
\newtheorem{thm}[dfn]{Theorem}
\newtheorem{rem}[dfn]{Remark}

\newtheorem{ex}[dfn]{Example}

\newcommand{\mf}{\mathfrak}
\DeclareMathOperator{\End}{\sf End}

\DeclareMathOperator{\Hom}{\sf Hom}

\DeclareMathOperator{\ch}{\sf ch}

\DeclareMathOperator{\id}{\sf id}

\DeclareMathOperator{\ad}{\sf ad}

\makeatletter

\@addtoreset{equation}{section}
\makeatother

\title{Kazama--Suzuki coset construction and its inverse}
\author{Ryo SATO}
\address{Institute of Mathematics, Academia Sinica, Taipei, Taiwan 10617}
\email{sato@gate.sinica.edu.tw}


\begin{document}

\maketitle

\begin{abstract}
We study the representation theory of
the Kazama--Suzuki coset vertex operator superalgebra
associated with the pair of a complex simple Lie algebra 
and its Cartan subalgebra.
In the case of type $A_{1}$,
B.L.\,Feigin, A.M.\,Semikhatov, and I.Yu.\,Tipunin
introduced another coset construction,
which is ``inverse'' of the Kazama--Suzuki coset construction.
In this paper we generalize the latter coset construction to arbitrary type
and establish a categorical equivalence between the categories of certain 
modules over an affine vertex operator algebra and the corresponding
Kazama--Suzuki coset vertex operator superalgebra.
Moreover, when the affine vertex operator algebra is regular,
we prove that the corresponding
Kazama--Suzuki coset vertex operator superalgebra
is also regular and the category of its ordinary modules
carries a braided monoidal category structure
by the theory of vertex tensor categories.
\end{abstract}

\setcounter{tocdepth}{1}
\tableofcontents

\section{Introduction}

Representation theory of vertex operator superalgebras plays a fundamental role
in the mathematical study of the corresponding $2$-dimensional
conformal field theories, see e.g.\,\cite{frenkel1992vertex}, 
\cite{nagatomo2005conformal}, \cite{huang2008vertex}.
When the symmetry of the Virasoro algebra on a given vertex operator superalgebra
is extended to that of $\mathbb{Z}/2\mathbb{Z}$-graded generalizations of the Virasoro algebra,
known as superconformal algebras (e.g.\,\cite{kac1988classification}),
the corresponding theory is naturally expected to provide
an example of superconformal field theories.
Though such generalizations are
originally motivated by purely physical applications (e.g.\,\cite{neveu1971factorizable}, \cite{ramond1971dual}),
the superconformal algebra symmetry
turns out to be non-trivially related with several areas of mathematics,
see e.g.\,\cite{kac2003quantum}, \cite{ben2008supersymmetry},
 \cite{eguchi2011notes}, \cite{witten2012notes}, \cite{kac2016representations}.

In this paper, we study the representation theory of a certain specific family 
of vertex operator superalgebras
with the $\mathcal{N}=2$ superconformal algebra symmetry,
which is given by the \emph{Kazama--Suzuki coset construction} \cite{KS89}
associated with the pair of a complex simple Lie algebra $\mf{g}$ and its Cartan subalgebra $\mf{h}$.
Our main tool to study such vertex operator superalgebras is
a generalization of the \emph{Feigin--Semikhatov--Tipunin coset construction}
introduced in \cite{FST98} for $\mf{g}=\mf{sl}_{2}$.
Roughly speaking, 
these two constructions are 
Heisenberg coset constructions (see \cite{creutzig2018schur} 
for the general theory and many examples), 
which are ``mutually inverse'' to each other:
\begin{enumerate}
\item the Kazama--Suzuki coset construction (``fermionization'')
$$(\text{\sf Affine VOAs})\to
(\mathcal{N}=2\ \text{\sf VOSAs}),$$
\item the Feigin--Semikhatov--Tipunin coset construction (``defermionization'')
$$(\mathcal{N}=2\ \text{\sf VOSAs})\to
(\text{\sf Affine VOAs}).$$
\end{enumerate}

One of the most important features of
the above bidirectional constructions is that
they relate the representation theory of $\mathbb{Z}/2\mathbb{Z}$-graded objects
to that of purely even\footnote{
We say that a $\mathbb{Z}/2\mathbb{Z}$-graded vector space
$V=V^{\bar{0}}\oplus V^{\bar{1}}$ is purely even if $V^{\bar{1}}=\{0\}$.
Throughout this paper, we regard all the vertex algebras as purely even vertex superalgebras.}
objects.
In fact, we establish 
a categorical equivalence between $\mathbb{C}$-linear additive categories of certain 
modules over affine vertex operator algebras and the corresponding $\mathcal{N}=2$ 
superconformal vertex operator superalgebras.

\subsection{Main result}

In order to explain the main result, we first give a brief review about 
the Kazama--Suzuki coset construction for the pair $(\mf{g},\mf{h})$ as above.

Let $k\in\mathbb{C}\setminus\{-h^{\vee}\}$, where $h^{\vee}$ is the dual Coxeter number
of $\mf{g}$.
We denote by $V^{k}(\mf{g})$ the universal affine vertex operator algebra (see Lemma \ref{UnivAff})
associated with the invariant bilinear form $kB$,
where $B$ is the normalized symmetric invariant bilinear form on $\mf{g}$.
We simply write $V_{\sf af}$ for $V^{k}(\mf{g})$ or its simple quotient $L_{k}(\mf{g})$.

As is described in \cite[\S2]{hosono1991lie},
the adjoint $\mf{h}$-actions on $\mf{g}$ and on 
the orthogonal complement $\mf{h}^{\perp}$ of $\mf{h}$ with respect to $B$
give rise to an injective vertex superalgebra homomorphism (see \S\ref{KS} for details)
\begin{equation}\label{reductive_action}
V^{k+h^{\vee}}(\mf{h})\to V_{\sf af}\otimes V^{+},
\end{equation}
where
$V^{k+h^{\vee}}(\mf{h})$ is the Heisenberg vertex operator algebra
associated with the bilinear form $(k+h^{\vee})B|_{\mf{h}\times\mf{h}}$ and $V^{+}$
is the Clifford vertex operator superalgebra associated with $\mf{h}^{\perp}$
and $B|_{\mf{h}^{\perp}\times\mf{h}^{\perp}}$.
In this paper, we write $\mathcal{H}^{+}$ for the image of (\ref{reductive_action})
and denote the corresponding commutant (usually referred to as \emph{coset}) 
vertex superalgebra by
\begin{equation*}
V_{\sf sc}:={\sf Com}\big(\mathcal{H}^{+},V_{\sf af}\otimes
V^{+}\big).
\end{equation*}
As a special case of \cite[Theorem 2.5]{hosono1991lie}
(see also \cite[\S3]{KS89}),
one can see that $V_{\sf sc}$
has the $\mathcal{N}=2$ superconformal algebra symmetry:
\begin{thm}[{\cite{hosono1991lie}}]
The coset vertex superalgebra $V_{\sf sc}$
carries a structure of 
an $\mathcal{N}=2$ superconformal vertex operator superalgebra\footnote{
It is clear that both $\mathcal{C}^{k}(\mf{g})$ and $\mathcal{C}_{k}(\mf{g})$ 
are of CFT type, that is,
every $L_{0}$-eigenvalue is non-negative
and the $L_{0}$-eigenspace of eigenvalue $0$ is $1$-dimensional.}
of central charge 
$$c_{\sf sc}:=\frac{k\,{\sf dim}\,\mf{g}}{k+h^{\vee}}
+\frac{1}{2}{\sf dim}\,\mf{h}^{\perp}-{\sf dim}\,\mf{h}$$
in the sense of \cite[Definition 1.1]{huang2002intertwining} (see also \cite[Definition 1.1]{Ad99}).
\end{thm}

Next we move on to the Feigin--Semikhatov--Tipunin coset construction.
In \cite{FST98}, B.L.\,Feigin, A.M.\,Semikhatov, and I.Yu.\,Tipunin introduced
a coset construction of the affine Lie algebra $\widehat{\mf{sl}}_{2}$
from the $\mathcal{N}=2$ superconformal algebra
and the lattice vertex superalgebra associated with the negative-definite lattice
$\sqrt{-1}\mathbb{Z}$.
After a while, D.\,Adamovi\'{c} reformulated their construction purely in terms 
of vertex superalgebras in \cite[\S4]{Ad99}.
The key result in this paper is an appropriate generalization of their construction
to arbitrary type as follows.

\begin{prp}[Proposition \ref{strong_FST}]
Assume that $k\in\mathbb{C}\setminus\{0,-h^{\vee}\}$.
Let $V^{-}$ be the lattice vertex superalgebra associated with
the negative-definite integral lattice $\sqrt{-1}\mathbb{Z}^{{\sf dim}\mf{h}}$.
Then, for a certain Heisenberg vertex subalgebra $\mathcal{H}^{-}$ of
$V_{\sf sc}\otimes V^{-}$ (see Lemma \ref{heis_minus} for the definition), 
there exists an explicit isomorphism
\begin{equation*}
f_{\sf FST}\colon V_{\sf af}\to
{\sf Com}\big(\mathcal{H}^{-},V_{\sf sc}\otimes V^{-}\big)
\end{equation*}
of purely even vertex operator superalgebras.
\end{prp}

For each coset $[\lambda]\in\mf{h}^{*}\!/Q$, where $Q$ is the root lattice of $\mf{g}$, 
we introduce certain 
$\mathbb{C}$-linear full subcategories 
$\mathscr{C}_{[\lambda]}(V_{\sf af})$
and $\mathscr{C}_{[\lambda_{\sf sc}]}(V_{\sf sc})$ of
weak $V_{\sf af}$-modules and weak $V_{\sf sc}$-modules, respectively
(see Definition \ref{harishchandra} and Proposition \ref{Funct}).
Now we can state our main result as follows:

\begin{thm}[Theorem \ref{MainThm}]
Let $k\in\mathbb{C}\setminus\{0,-h^{\vee}\}$ and $\lambda\in\mf{h}^{*}$.
Then the following $\mathbb{C}$-linear functors
\begin{align*}
\Omega^{+}(\lambda)\colon
\mathscr{C}_{[\lambda]}(V_{\sf af})
\to\mathscr{C}_{[\lambda_{\sf sc}]}(V_{\sf sc}),\ \ 
\Omega^{-}(\lambda_{\sf sc})\colon
\mathscr{C}_{[\lambda_{\sf sc}]}(V_{\sf sc})
\to\mathscr{C}_{[\lambda]}(V_{\sf af})
\end{align*}
defined in \S\ref{Omega_def} are mutually quasi-inverse to each other.
Moreover, these functors satisfy 
a certain equivariance property with respect to the \emph{spectral flow automorphisms}
(see \S\ref{SpecFlowEquiv} for details).
\end{thm}

We should mention that the above result for $\mf{g}=\mf{sl}_{2}$ is 
essentially the same as
\cite[Theorem 4.4, 7.7, and Corollary 6.3]{sato2016equivalences} and the original idea in this case
goes back to \cite[Theorem IV.\,10]{FST98}.
It is worth noting that, when $\mf{g}\neq\mf{sl}_{2}$, 
the same proof as \cite[Theorem 4.4]{sato2016equivalences}
does not work due to
the existence of a non-trivial lattice vertex subsuperalgebra of $V_{\sf sc}$
(see \S\ref{lattice_K}).
See Remark \ref{difficulty} for details.

\subsection{Regular cases}

It is well-known that
the simple affine vertex operator algebra
$L_{k}(\mf{g})$ is regular if $k$ is a positive integer (see \cite[Theorem 3.7]{dong1997regularity}).
In this case, we obtain further information about
the simple vertex operator superalgebra
$$\mathcal{C}_{k}(\mf{g})
:={\sf Com}\big(\mathcal{H}^{+},L_{k}(\mf{g})\otimes V^{+}\big).$$

We first determine the bicommutant of the Heisenberg vertex subalgebra
$\mathcal{H}^{+}$ and $\mathcal{H}^{-}$
 in $L_{k}(\mf{g})\otimes V^{+}$ and $\mathcal{C}_{k}(\mf{g})\otimes V^{-}$, respectively.

\begin{prp}[Proposition \ref{plus_extension} and \ref{minus_extension}]
Assume that $k$ is a positive integer.
\begin{enumerate}
\item The bicommutant vertex superalgebra
$$\mathcal{E}^{+}:={\sf Com}\Big({\sf Com}\big(\mathcal{H}^{+},L_{k}(\mf{g})\otimes V^{+}\big),
L_{k}(\mf{g})\otimes V^{+}\Big)$$
is purely even and isomorphic to the lattice vertex algebra
associated with $$\sqrt{k+h^{\vee}}Q_{l},$$ where $Q_{l}$ is the even integral sublattice of 
the root lattice $Q$ spanned by long roots (cf.\,\cite[Proposition 4.11]{dong2011parafermion}).
\item The bicommutant vertex superalgebra
$$\mathcal{E}^{-}:={\sf Com}\Big(
{\sf Com}\big(\mathcal{H}^{-},\mathcal{C}_{k}(\mf{g})\otimes V^{-}\big),
\mathcal{C}_{k}(\mf{g})\otimes V^{-}\Big)$$
is isomorphic to 
the lattice vertex superalgebra associated with 
$$\sqrt{-(k+h^{\vee})}Q_{l}\oplus\mathbb{Z}^{\frac{1}{2}{\sf dim}\mf{h}^{\perp}-{\sf dim}\mf{h}}.$$
\end{enumerate}
\end{prp}

Next, by using \cite[Corollary 2]{miyamoto2015c_2}
and \cite[Theorem 4.12]{creutzig2018schur}
(cf.\,\cite[Theorem 5.24]{carnahan2016regularity}),
we prove the following.

\begin{thm}[Theorem \ref{regular}]\label{Reg}
When $k$ is a positive integer, both 
$\mathcal{C}_{k}(\mf{g})$ and its even part $\mathcal{C}_{k}(\mf{g})^{\bar{0}}$ are regular.
\end{thm}

At last, we obtain the following theorem as a special case of 
\cite[Theorem 3.65]{creutzig2017tensor},
which is based on the theory of vertex tensor categories developed by
Y.-Z.\,Huang and J.\,Lepowsky (see 
\cite{huang1995theory1},
\cite{huang1995theory2},
\cite{huang1995theory3},
\cite{huang1995theory4},
\cite{huang2005differential}, and references therein).

\begin{thm}[Theorem \ref{BTC}]\label{Braided}
When $k$ is a positive integer, the semisimple $\mathbb{C}$-linear abelian
category of $\mathbb{Z}/2\mathbb{Z}$-graded
ordinary $\mathcal{C}_{k}(\mf{g})$-modules
carries a braided monoidal category structure\footnote{
Note that the symmetric center (see \cite[Definition 2.9]{muger2003structure}) 
of this braided monoidal category
contains at least two non-isomorphic simple objects,
the monoidal unit $\mathcal{C}_{k}(\mf{g})$ and
its $\mathbb{Z}/2\mathbb{Z}$-parity reversed object.
In particular, it is not modular in the sense of \cite[Definition 3.1]{muger2003structure}.}
induced by \cite[Theorem 3.7]{huang2005differential} 
and \cite[Theorem 3.65]{creutzig2017tensor}.
\end{thm}

We should note that Theorem \ref{Reg} (resp.\,Theorem \ref{Braided})
for $\mf{g}=\mf{sl}_{2}$ is proved by 
D.\,Adamovi\'{c} in \cite[Theorem 8.1]{Ad01}
(resp.\,by Y.-Z.\,Huang and
A.\,Milas in \cite[Theorem 4.8]{huang2002intertwining}).

\subsection{Structure of the paper}

We organize this paper as follows.
In \S2 we define the coset vertex operator superalgebra $V_{\sf sc}$
and its ``Cartan subalgebra'' $\mf{h}_{\sf sc}$ in an explicit way.
Our main tool, a generalization of the Feigin--Semikhatov--Tipunin
coset construction, is introduced in \S\ref{Defermi}.
In \S3, after introducing an appropriate category of modules
and two families of functors $\Omega^{+}$ and $\Omega^{-}$,
we give the precise statement of our main result.
The proof of the main theorem is given in \S4.
At last, in \S5, we discuss the case that the level $k$ is a positive integer.

\medskip
{\bf Acknowledgments: }
The author would like to express his gratitude to 
Shun-Jen Cheng, Ching Hung Lam, Masahiko Miyamoto, Hiromichi Yamada,
and Hiroshi Yamauchi for helpful discussions and valuable advice.
He also would like to express his appreciation to 
Tomoyuki Arakawa, Kenichiro Tanabe, Shintarou Yanagida, and many others
for constructive comments.

\section{Fermionization and Defermionization}

\subsection{Kazama--Suzuki coset construction}\label{KS}

In this subsection, we review the vertex superalgebraic formulation
of the Kazama-Suzuki coset constructions.
See \cite{KS89} and \cite{hosono1991lie} for details.

Let $\mf{g}$ be a simple complex Lie algebra 
with a fixed Borel subalgebra $\mf{b}$.
Let $\mf{h}$ be the corresponding Cartan subalgebra
and $\Delta$ (resp.\,$\Delta^{+}, \Pi$) the set of (resp.\,positive, simple) roots 
associated with $(\mf{g},\mf{b})$.
We write $B=B(?,?)\colon\mf{g}\times\mf{g}\to\mathbb{C}$ for the normalized 
symmetric invariant bilinear form on $\mf{g}$.

Define a positive-definite integral lattice $L^{+}$ of rank $N:=|\Delta^{+}|$ by
$$L^{+}:=\bigoplus_{\alpha\in\Delta^{+}}\mathbb{Z}\alpha^{+},\ \ 
\langle\alpha^{+},\beta^{+}\big\rangle:=\delta_{\alpha,\beta}\ (\alpha,\beta\in\Delta^{+})$$
and take a bimultiplicative map $\epsilon^{+}\colon L^{+}\times L^{+}\to\{\pm1\}$
in the same way as $\epsilon\colon\mathbb{Z}^{N}\times\mathbb{Z}^{N}\to\{\pm1\}$
in \cite[Example 5.5a]{kac97V}\footnote{
Our discussion below does not depend on the choice 
of the lattice isomorphism $L^{+}\simeq\mathbb{Z}^{N}$.}.
Then we write $V^{+}$ for the lattice vertex superalgebra
associated with $(L^{+},\epsilon^{+})$, which is
isomorphic to the $N$-times tensor product of 
the charged free fermions.
Unless otherwise specified,
we follow the standard notation of lattice vertex superalgebras 
used in \cite{kac97V}.

Let $h^{\vee}$ be the dual Coxeter number of $\mf{g}$ and 
$k\in\mathbb{C}\setminus\{-h^{\vee}\}$. 
We denote by $V_{\sf af}$ the universal affine vertex operator algebra $V^{k}(\mf{g})$
of level $k$ or its simple quotient vertex operator algebra $L_{k}(\mf{g})$.
Based on the setting of the Kazama--Suzuki coset construction,
we consider a Heisenberg vertex subalgebra of $V_{\sf af}\otimes V^{+}$ as follows.

\begin{lem}
For $\alpha\in\Delta$, we set
$$H_{\alpha}:=\frac{(\alpha,\alpha)}{2}\alpha^{\vee}\in\mf{h},$$
where $(?,?)\colon\mf{h}^{*}\times\mf{h}^{*}\to\mathbb{C}$ is 
the normalized bilinear form induced by 
the restriction of $B$ to $\mf{h}\times\mf{h}$ and
$\alpha^{\vee}\in\mf{h}$ is the coroot of $\alpha$.
We also set
\begin{align*}
H^{+}_{\alpha}
&:=H_{\alpha,-1}\mathbf{1}_{\sf af}\otimes\mathbf{1}_{V^{+}}
-\mathbf{1}_{\sf af}\otimes\,\alpha_{\sf f,-1}\mathbf{1}_{V^{+}}\in (V_{\sf af}\otimes V^{+})^{\bar{0}},
\end{align*}
where $\mathbf{1}_{\sf af}$ is the vacuum vector of $V_{\sf af}$ and
$$\alpha_{\sf f}:=\sum_{\beta\in\Delta^{+}}
(\alpha,\beta)\beta^{+}\in L^{+}\otimes_{\mathbb{Z}}\mathbb{Q}.$$
Then the vertex subalgebra $\mathcal{H}^{+}$ of
$V_{\sf af}\otimes V^{+}$ generated by the set
$\{H^{+}_{\alpha}\,|\,\alpha\in \Pi\}$
is the Heisenberg vertex algebra of rank $\ell:={\sf dim}\,\mf{h}$.
\end{lem}

\begin{proof}
Since we have
$$\langle \alpha_{\sf f},\beta_{\sf f}\rangle
=\sum_{\gamma\in\Delta^{+}}
(\alpha,\gamma)(\gamma,\beta)
=h^{\vee}(\alpha,\beta)$$
for $\alpha,\beta\in\Pi$,
there exists a unique vertex algebra isomorphism
$V^{k+h^{\vee}}(\mf{h})\to\mathcal{H}^{+}$
such that
$H_{\alpha,-1}\mathbf{1}_{\sf af}\mapsto H^{+}_{\alpha}$
for $\alpha\in \Pi$.
\end{proof}

Let $\omega^{\sf af}$, $\omega^{\sf f}$, and $\omega^{\sf b}$ be
the standard conformal vectors of central charge
$c_{\sf af}:=\frac{k\,{\sf dim}\,\mf{g}}{k+h^{\vee}}$,
$N$, and $\ell$ of $V_{\sf af}$, $V^{+}$, and $\mathcal{H}^{+}$, respectively 
(see e.g.\,\cite{kac97V} for details).

\begin{lem}
The coset vertex superalgebra 
$$V_{\sf sc}:={\sf Com}
\bigl(\mathcal{H}^{+}\!,V_{\sf af}\otimes V^{+}\bigr)$$ together with
$\omega^{\sf sc}:=\omega^{\sf af}\otimes\mathbf{1}_{V^{+}}
+\mathbf{1}_{\sf af}\otimes \omega^{\sf f}-\omega^{\sf b}$
forms a vertex operator superalgebra
of central charge 
$c_{\sf sc}:=c_{\sf af}+N-\ell.$
\end{lem}

\begin{proof}
It follows from the general theory of coset vertex superalgebras (see e.g.\,\cite[Theorem 5.1]{frenkel1992vertex}).
\end{proof}

\begin{prp}\label{simplequot}
The coset vertex operator superalgebra
$$\mathcal{C}_{k}(\mf{g}):={\sf Com}
\bigl(\mathcal{H}^{+}\!,L_{k}(\mf{g})\otimes V^{+}\bigr)$$
is a unique simple quotient vertex operator superalgebra of 
$$\mathcal{C}^{k}(\mf{g}):={\sf Com}
\bigl(\mathcal{H}^{+}\!,V^{k}(\mf{g})\otimes V^{+}\bigr).$$
\end{prp}

\begin{proof}
Let $I$ be the maximal proper ideal of $V^{k}(\mf{g})$.
Since $V^{k}(\mf{g})\otimes V^{+}$ and $I\otimes V^{+}$ are completely reducible
as weak $\mathcal{H}^{+}$-modules by the same argument as in \cite[\S1.7]{frenkel1989vertex},
it follows that $\mathcal{C}_{k}(\mf{g})$ is a quotient of $\mathcal{C}^{k}(\mf{g})$.
The simplicity and uniqueness follow from 
the same argument of \cite[Proposition 3.2]{creutzig2018schur}
(see also \cite[Lemma 2.1]{arakawa2017orbifolds})
and the fact that $\mathcal{C}^{k}(\mf{g})$ is of CFT type, respectively.
\end{proof}

\begin{lem}\label{coset1}
The following homomorphism
\begin{equation*}
\iota_{+}\colon V_{\sf sc}\otimes\mathcal{H}^{+}
\to V_{\sf af}\otimes V^{+};\,A\otimes B\mapsto A_{-1}B
\end{equation*}
of vertex operator superalgebras is injective.
\end{lem}

\begin{proof}
Since $\mathcal{H}^{+}$ is simple and
$V_{\sf af}\otimes V^{+}$ is completely reducible as
a weak $\mathcal{H}^{+}$-module,
the injectivity of $\iota_{+}$ follows.
\end{proof}

\subsection{Cartan subalgebra}

In this subsection, we introduce a Heisenberg vertex subalgebra of $V_{\sf sc}$.
From now on, we always assume that $k\in\mathbb{C}\setminus\{0,-h^{\vee}\}$.

\begin{lem}\label{cosetgen}
The following set
\begin{align*}
\left\{\left.J_{\alpha}:=
\frac{1}{k}H_{\alpha,-1}\mathbf{1}_{\sf af}\otimes\mathbf{1}_{V^{+}}
+\mathbf{1}_{\sf af}\otimes\alpha^{+}_{-1}\mathbf{1}_{V^{+}}\in(V_{\sf af}\otimes V^{+})^{\bar{0}}
\,\right|\,
\alpha\in\Delta^{+}\right\}
\end{align*}
generates a vertex subalgebra $\mathcal{H}_{\sf sc}$
of $V_{\sf sc}$,
which is isomorphic to the Heisenberg vertex algebra of rank $N$. 
\end{lem}

\begin{proof}
Since it is easy to verify that each
$J_{\alpha}$ lies in $V_{\sf sc}$,
we only need to check the non-degeneracy.
By a direct computation,
we have
\begin{align*}
J_{\alpha}(z)J_{\beta}(w)\sim
\frac{g_{\alpha,\beta}}{(z-w)^{2}}:=
\frac{\frac{(\alpha,\beta)}{k}+\delta_{\alpha,\beta}}{(z-w)^{2}},\ 
\end{align*}
for any $\alpha,\beta\in\Delta^{+}$.
By using the formula
\begin{equation*}
\sum_{\alpha\in\Delta^{+}}(\lambda,\alpha)\alpha
=h^{\vee}\lambda
\end{equation*}
for $\lambda\in\mf{h}^{*}$, we can verify that
$$(g^{*}_{\alpha,\beta})_{\alpha,\beta\in\Delta^{+}}
:=\left(-\frac{(\alpha,\beta)}
{k+h^{\vee}}+\delta_{\alpha,\beta}\right)_{\alpha,\beta\in\Delta^{+}}$$
is the inverse matrix of $(g_{\alpha,\beta})_{\alpha,\beta\in\Delta^{+}}.$
This completes the proof.
\end{proof}

By the inverse matrix in the above proof, we obtain the following.

\begin{cor}\label{Jstar_OPE}
We set
\begin{align*}
J_{\alpha}^{*}:
=\sum_{\beta\in\Delta^{+}}g^{*}_{\alpha,\beta}J_{\beta}
=\frac{1}{k+h^{\vee}}H^{+}_{\alpha}
+\mathbf{1}_{\sf af}\otimes\,\alpha_{-1}^{+}\mathbf{1}_{V^{+}}\in\mathcal{H}_{\sf sc}
\end{align*}
for $\alpha\in\Delta^{+}$.
Then we have
$$J_{\alpha}^{*}(z)J_{\beta}(w)\sim
\frac{\delta_{\alpha,\beta}}{(z-w)^{2}},\ \ 
J_{\alpha}^{*}(z)J^{*}_{\beta}(w)\sim
\frac{g^{*}_{\alpha,\beta}}{(z-w)^{2}}$$
for any $\alpha,\beta\in\Delta^{+}$.
\end{cor}

\begin{lem}
We set
$$P^{\vee}_{\sf sc}:=\bigoplus_{\alpha\in\Delta^{+}}\mathbb{Z}J_{\alpha}^{*},\ \ 
\mf{h}_{\sf sc}:=\bigoplus_{\alpha\in\Delta^{+}}\mathbb{C}J_{\alpha}^{*},$$
and
$Q_{\sf sc}:=\big\{\lambda\in(\mf{h}_{\sf sc})^{*}\,
\big|\,\lambda(H)\in\mathbb{Z}
\text{ for any }H\in P^{\vee}_{\sf sc}\big\}.$
Then $V_{\sf sc}$ decomposes into a direct sum
$$V_{\sf sc}=\bigoplus_{\gamma\in Q_{\sf sc}}(V_{\sf sc})_{\gamma}$$
of $\mf{h}_{\sf sc}$-weight subspaces, where
$$(V_{\sf sc})_{\gamma}:=\{A\in V_{\sf sc}\,|\,H_{0}A
=\gamma(H)A\text{ for any }H\in\mf{h}_{\sf sc}\}.$$
\end{lem}

\begin{proof}
It follows from the fact that the operator $J_{\alpha,0}^{*}$
coincides with $\mathbf{1}_{\sf af}\otimes\,\alpha_{0}^{+}$ on $V_{\sf sc}$
for any $\alpha\in\Delta^{+}$.
\end{proof}

\subsection{Feigin--Semikhatov--Tipunin coset construction}\label{Defermi}

In \cite{FST98}, B.L.\,Feigin, A.M.\,Semikhatov, and I.Yu.\,Tipunin
studied a coset construction
of $V^{k}(\mf{sl}_{2})$ in the tensor product of $\mathcal{C}^{k}(\mf{sl}_{2})$
and the lattice vertex superalgebra associated with $\sqrt{-1}\mathbb{Z}$
(see also \cite[\S5]{Ad99}).
In this subsection, we generalize their coset construction to arbitrary type.

We define a negative-definite integral lattice $L^{-}$ of rank $\ell$ by
$$L^{-}:=\bigoplus_{\alpha\in \Pi}\mathbb{Z}\alpha^{-},\ \ 
\langle\alpha^{-},\beta^{-}\rangle:=-\delta_{\alpha,\beta}$$
and a bimultiplicative map $\epsilon^{-}\colon L^{-}\times L^{-}\to\{\pm1\}$ by
$$\epsilon^{-}(\alpha^{-},\beta^{-}):=
-\epsilon^{+}(\alpha^{+},\beta^{+})$$
for $\alpha,\beta\in\Pi$.
We write $V^{-}$ for the lattice vertex superalgebra associated with $(L^{-},\epsilon^{-})$\footnote{
One can easily verify that $\epsilon^{-}$ satisfies the condition 
$$\epsilon^{-}(\zeta,\zeta')\epsilon^{-}(\zeta',\zeta)
=(-1)^{\langle\zeta,\zeta'\rangle+\langle\zeta,\zeta\rangle\langle\zeta',\zeta'\rangle}$$
for any $\zeta,\zeta'\in L^{-}$ in \cite[Theorem 5.5 (b)]{kac97V}.}.
Since the lattice $L^{-}$ is negative-definite,
the vertex superalgebra $V^{-}$ has a natural conformal structure of
central charge $\ell$ and the set of the $L_{0}$-eigenvalues is not 
bounded below.

\begin{lem}\label{heis_minus}
For each $\alpha\in\Delta^{+}$, we set
$$H^{-}_{\alpha}:=
\begin{cases}
J^{*}_{\alpha}\otimes\mathbf{1}_{V^{-}}
+\mathbf{1}_{\sf sc}\otimes\,\alpha^{-}_{-1}\mathbf{1}_{V^{-}} &
\text{ if }\alpha\in\Pi,\\
J^{*}_{\alpha}\otimes\mathbf{1}_{V^{-}} &
\text{ if }\alpha\notin\Pi,
\end{cases}$$
where $\mathbf{1}_{\sf sc}$ is the vacuum vector of $V_{\sf sc}$.
Then the vertex subalgebra $\mathcal{H}^{-}$ of $V_{\sf sc}\otimes V^{-}$
generated by $\{H^{-}_{\alpha}\,|\,\alpha\in\Delta^{+}\}$
is the Heisenberg vertex algebra of rank $N$.
\end{lem}

\begin{proof}
By Corollary \ref{Jstar_OPE}, we have
$$H^{-}_{\alpha}(z)H^{-}_{\beta}(w)
\sim
\frac{G_{\alpha,\beta}}{(z-w)^{2}}:=
\begin{cases}
\displaystyle\frac{g^{*}_{\alpha,\beta}-\delta_{\alpha,\beta}}{(z-w)^{2}}
&\text{ if }\alpha\in\Pi\text{ or }\beta\in\Pi,\\
\displaystyle
\frac{g^{*}_{\alpha,\beta}}{(z-w)^{2}} 
&\text{ if }\alpha,\beta\notin\Pi
\end{cases}$$
for $\alpha,\beta\in\Delta^{+}$.
We define an $N\times N$ matrix $(G^{*}_{\alpha,\beta})_{\alpha,\beta\in\Delta^{+}}$ by
\begin{equation}\label{dualG}
G^{*}_{\alpha,\beta}:=
\begin{cases}
-kc_{\alpha,\beta}-\delta_{\alpha,\beta} & \text{ if }\alpha,\beta\in\Pi,\\
-\beta(\alpha^{*}) & \text{ if }\alpha\in\Pi\text{ and }\beta\notin\Pi,\\
-\alpha(\beta^{*}) & \text{ if }\alpha\notin\Pi\text{ and }\beta\in\Pi,\\
\delta_{\alpha,\beta} & \text{ if }\alpha,\beta\notin\Pi,
\end{cases}
\end{equation}
where $(c_{\alpha,\beta})_{\alpha,\beta\in\Pi}$ is the inverse
of the symmetrized Cartan matrix $\big((\alpha,\beta)
\big)_{\alpha,\beta\in\Pi}$ and
$\alpha^{*}\in\mf{h}$ is the fundamental coweight
defined by $\beta(\alpha^{*})=\delta_{\alpha,\beta}$ for $\beta\in\Pi$.
Then we can verify that 
$(G^{*}_{\alpha,\beta})_{\alpha,\beta\in\Delta^{+}}$ is the inverse matrix 
of $(G_{\alpha,\beta})_{\alpha,\beta\in\Delta^{+}}$ and thus the latter is non-singular.
\end{proof}

Let $Q$ be the root lattice\footnote{
The root lattice $Q$ with the normalized bilinear form $(?,?)$
is integral (resp.\,even) if and only if $\mf{g}$ is not of type $G_{2}$ (resp.\,is of type $A,D,E$).}
of $\mf{g}$.
We define $\mathbb{Z}$-linear maps $f^{\pm}_{\sf af}\colon Q\to L^{\pm}$ by
\begin{equation*}
f^{\pm}_{\sf af}(\gamma):=
\sum_{\alpha\in\Pi}\gamma(\alpha^{*})\alpha^{\pm}
\end{equation*}
for $\gamma\in Q$.
Then we obtain the next lemma as a generalization 
of the ``anti-Kazama--Suzuki mapping'' for $\mf{g}=\mf{sl}_{2}$
in \cite[Lemma III.5]{FST98} and \cite[Lemma 5.1]{Ad99} to arbitrary type.

\begin{lem}\label{FSTmap}
There exists a unique vertex superalgebra homomorphism
\begin{equation}\label{defermionize}
f_{\sf FST}\colon V^{k}(\mf{g})\to
{\sf Com}\big(\mathcal{H}^{-},\mathcal{C}^{k}(\mf{g})\otimes V^{-}\big).
\end{equation}
such that 
\begin{align*}
X_{\alpha,-1}\mathbf{1}_{\sf af}&\mapsto 
\widetilde{X}_{\alpha}:=X_{\alpha,-1}\mathbf{1}_{\sf af}\otimes
\mathrm{e}^{f^{+}_{\sf af}(\alpha)}
\otimes\mathrm{e}^{f^{-}_{\sf af}(\alpha)},\\
H_{\beta,-1}\mathbf{1}_{\sf af}&\mapsto
\widetilde{H}_{\beta}
:=k\big(J_{\beta}\otimes\mathbf{1}_{V^{-}}
+\mathbf{1}_{\sf sc}\otimes\beta^{-}_{-1}\mathbf{1}_{V^{-}}\big),
\end{align*}
for any $\alpha\in\Delta$ and $\beta\in\Pi$,
where
$X_{\alpha}\in\mf{g}$ is a root vector associated with
the root $\alpha$ and normalized by $[X_{\alpha},X_{-\alpha}]=\alpha^{\vee}$.
\end{lem}

\begin{proof}
First we prove the above assignment
induces a vertex superalgebra homomorphism
from $V^{k}(\mf{g})$ to $V^{k}(\mf{g})\otimes V^{+}\otimes V^{-}$.
It suffices to show that the above elements in the right-hand side obey the same OPE as
in the left-hand side.
By the general theory of lattice vertex superalgebras,
we have
\begin{align}\label{VO_OPE}
&\mathrm{e}^{\xi}
_{n}\mathrm{e}^{\xi'}
=\begin{cases}
0 & \text{ if }n\geq-\langle\xi,\xi'\rangle,\\
\epsilon(\xi,\xi')\mathrm{e}^{\xi+\xi'} & \text{ if }n=-\langle\xi,\xi'\rangle-1,\\
\epsilon(\xi,\xi')\xi_{-1}\mathrm{e}^{\xi+\xi'} & \text{ if }n=-\langle\xi,\xi'\rangle-2
\end{cases}
\end{align}
for elements $\xi,\xi'$ of a general non-degenerate integral lattice,
where $\epsilon$ is an appropriate $2$-cocycle of the lattice
(see e.g.\,\cite[\S5.5]{kac97V}).
In particular, since we have 
$\langle f^{+}_{\sf af}(\gamma),f^{+}_{\sf af}(\gamma')\rangle
+\langle f^{-}_{\sf af}(\gamma),f^{-}_{\sf af}(\gamma')\rangle=0$
for any $\gamma,\gamma'\in Q$,
we obtain
\begin{equation}\label{Triv_OPE}
\left(\mathrm{e}^{f^{+}_{\sf af}(\gamma)}(z)
\otimes\mathrm{e}^{f^{-}_{\sf af}(\gamma)}(z)\right)
\left(\mathrm{e}^{f^{+}_{\sf af}(\gamma')}(w)
\otimes\mathrm{e}^{f^{-}_{\sf af}(\gamma')}(w)\right)\sim0.
\end{equation}
Then, by (\ref{VO_OPE}), (\ref{Triv_OPE}), and some computations, we obtain
\begin{align*}
\widetilde{X}_{\alpha}(z)
\widetilde{X}_{\beta}(w)
&\sim
\frac{\widetilde{[X_{\alpha},X_{\beta}]}(w)}{z-w}
+\frac{\delta_{\alpha+\beta,0}k\id\otimes\id\otimes\id}{(z-w)^{2}}
\end{align*}
for any $\alpha,\beta\in\Delta$.
In addition, we also obtain
$$\widetilde{H}_{\alpha}(z)\widetilde{X}_{\beta}(w)
\sim\frac{\widetilde{[H_{\alpha},X_{\beta}]}(w)}{z-w},\ 
\widetilde{H}_{\alpha}(z)\widetilde{H}_{\beta}(w)\sim
\frac{k(\alpha,\beta)\id\otimes\id\otimes\id}{(z-w)^{2}}.$$
Thus we get a homomorphism
from $V^{k}(\mf{g})$ to $V^{k}(\mf{g})\otimes V^{+}\otimes V^{-}$.

Next we need to prove that the image of the homomorphism
is actually contained in the commutant of
$\mathcal{H}^{-}$ in $\mathcal{C}^{k}(\mf{g})\otimes V^{-}$.
It suffices to check that
$(H_{\alpha}^{+}(z)\otimes\id_{V^{-}})\widetilde{X}_{\beta}(w)\sim0$
and $H_{\alpha}^{-}(z)\widetilde{X}_{\beta}(w)\sim0$
for any $\alpha\in\Delta^{+}$ and $\beta\in\Delta$.
Since they are verified by straightforward computations, we omit the detail.
\end{proof}

\begin{cor}
The homomorphism (\ref{defermionize}) gives rise to conformal vertex superalgebra\footnote{
Note that $V^{-}$ is not a vertex operator superalgebra.} homomorphisms
$\iota_{-}\colon V^{k}(\mf{g})\otimes\mathcal{H}^{-}\to 
\mathcal{C}^{k}(\mf{g})\otimes V^{-}$
and
$\iota_{\sf af}\colon V^{k}(\mf{g})\otimes
\mathcal{H}^{+}\otimes\mathcal{H}^{-}\to 
V^{k}(\mf{g})\otimes V^{+}\otimes V^{-}.$
\end{cor}

As a matter of fact, the following stronger statement holds.

\begin{prp}\label{strong_FST}
The homomorphism (\ref{defermionize}) is an isomorphism.
In particular, it descends to their simple quotients
and gives rise to 
conformal vertex superalgebra homomorphisms
$\iota_{-}\colon L_{k}(\mf{g})\otimes\mathcal{H}^{-}\to 
\mathcal{C}_{k}(\mf{g})\otimes V^{-}$
and $\iota_{\sf af}\colon L_{k}(\mf{g})\otimes
\mathcal{H}^{+}\otimes\mathcal{H}^{-}\to 
L_{k}(\mf{g})\otimes V^{+}\otimes V^{-}.$
\end{prp}

The proof of Proposition \ref{strong_FST} is postponed to Corollary \ref{FSTisom}.

\section{Main Theorem}

In this section, by using the notation defined in the previous section,
we give the precise statement of our main result (Theorem \ref{MainThm}).
Recall that the pair $(V_{\sf af},V_{\sf sc})$ stands for 
$\big(V^{k}(\mf{g}),\mathcal{C}^{k}(\mf{g})\big)$
or $\big(L_{k}(\mf{g}),\mathcal{C}_{k}(\mf{g})\big)$.
Throughout this section, we assume that $k\in\mathbb{C}\setminus\{0,-h^{\vee}\}$
and the letter $V$ always stands for 
$V_{\sf af}$ or $V_{\sf sc}$.

\subsection{Setting of the category of modules}

For $V=V_{\sf af}$ \big(resp.\,$V_{\sf sc}$\big),
we set
\begin{align*}
&P^{\vee}_{V}:=P^{\vee}
\ (\text{resp.\,}P^{\vee}_{\sf sc}),\ \ 
\mf{t}_{V}:=\mf{h}\ (\text{resp.\,}\mf{h}_{\sf sc}),\ \ 
Q_{V}:=Q\ (\text{resp.\,}Q_{\sf sc}),
\end{align*}
where $P^{\vee}$ is the coweight lattice in $\mf{h}$.

\begin{dfn}\label{harishchandra}
Let $\lambda\in\mf{t}_{V}^{*}$ and set
$$[\lambda]:=\{\lambda\}+Q_{V}\in\mf{t}_{V}^{*}/Q_{V}.$$
Then we define $\mathscr{C}_{[\lambda]}(V)$ to be a full subcategory of
$V\text{-\sf gMod}$ 
(see \S\ref{weakmod} for the definition) whose object
$M=M^{\bar{0}}\oplus M^{\bar{1}}$ satisfies the following conditions:
\begin{enumerate}
\item For each $\bar{i}\in\mathbb{Z}/2\mathbb{Z}$,
the $\mathbb{Z}/2\mathbb{Z}$-homogeneous subspace
$M^{\bar{i}}$ decomposes into a direct sum
$$M^{\bar{i}}=\bigoplus_{\mu\in[\lambda]}M_{\mu}^{\bar{i}}$$
of $\mf{t}_{V}$-weight spaces, where
$$M_{\mu}^{\bar{i}}:=\{v\in M^{\bar{i}}\,|\,H_{0}v=\mu(H)v\text{ for any }H\in\mf{t}_{V}\}.$$
\item For each $\bar{i}\in\mathbb{Z}/2\mathbb{Z}$ and $\mu\in[\lambda]$, 
the $\mf{t}_{V}$-weight space $M_{\mu}^{\bar{i}}$
further decomposes into a direct sum
$$M_{\mu}^{\bar{i}}=\bigoplus_{\Delta\in\mathbb{C}}M_{\mu}^{\bar{i}}(h)$$
of finite-dimensional subspaces, where
$$M_{\mu}^{\bar{i}}(h):=\{v\in M_{\mu}^{\bar{i}}\,|\,
(L_{0}-h)^{n}v=0\text{ if }n\gg0\}.$$
In addition, for any $h\in\mathbb{C}$,
we have $M_{\mu}^{\bar{i}}(h-r)=\{0\}$ if $r\gg0$.
\end{enumerate}
\end{dfn}

We write $M_{\mu}$
and $M_{\mu}(h)$
for $M_{\mu}^{\bar{0}}\oplus M_{\mu}^{\bar{1}}$
and $M_{\mu}^{\bar{0}}(h)\oplus M_{\mu}^{\bar{1}}(h)$, respectively.

\begin{dfn}
For an object $M$ of $\mathscr{C}_{[\lambda]}(V)$,
we define the \emph{string function of $M$ through $\mu\in[\lambda]$} by
$$s^{\mu}_{M}(q):=\sum_{\Delta\in\mathbb{C}}
{\sf dim}_{\mathbb{C}}\big(M_{\mu}(\Delta)\big)q^{\Delta-\frac{1}{24}c}$$
and the \emph{formal character of $M$} by
$$\ch(M):=\sum_{\mu\in[\lambda]}s^{\mu}_{M}(q)\mathrm{e}^{\mu}.$$
\end{dfn}

\subsection{Definition of the functors $\Omega^{+}$ and $\Omega^{-}$}\label{Omega_def}

In this subsection, we assume Proposition \ref{strong_FST} 
whenever we consider the case of $V=\mathcal{C}_{k}(\mf{g})$.
Note that the other cases are proved independently of Proposition \ref{strong_FST}.

We set 
$$\mf{t}^{+}:=\bigoplus_{\alpha\in\Pi}\mathbb{C}H_{\alpha}^{+}
\subsetneq\mathcal{H}^{+},\ \ 
\mf{t}^{-}:=\bigoplus_{\beta\in\Delta^{+}}\mathbb{C}H_{\beta}^{-}
\subsetneq\mathcal{H}^{-}.$$
From now on, we fix the following linear isomorphisms
\begin{align*}
&\mf{t}^{+}\xrightarrow{\simeq}\mf{t}_{V_{\sf af}}=\mf{h};\,
H_{\alpha}^{+}\mapsto H_{\alpha},\ \ 
\mf{t}^{-}\xrightarrow{\simeq}\mf{t}_{V_{\sf sc}}=\mf{h}_{\sf sc};\,
H_{\beta}^{-}\mapsto 
J_{\beta}^{*}
\end{align*}
and denote the induced linear isomorphisms by
\begin{align*}
\nu^{+}\colon\mf{t}_{V_{\sf af}}^{*}=\mf{h}^{*}\xrightarrow{\simeq}(\mf{t}^{+})^{*},\ \ 
\nu^{-}\colon\mf{t}_{V_{\sf sc}}^{*}=(\mf{h}_{\sf sc})^{*}\xrightarrow{\simeq}(\mf{t}^{-})^{*}.
\end{align*}

\begin{dfn}\label{Omegafunctor}
Let $\lambda,\mu\in\mf{t}_{V}^{*}$ and $M$ be an object of $\mathscr{C}_{[\lambda]}(V)$.
\begin{enumerate}
\item When $V=V_{\sf af}$, we define
a $\mathbb{Z}/2\mathbb{Z}$-graded weak $V_{\sf sc}$-module by
\begin{align*}
\Omega^{+}(\mu)(M):=\{v\in M\otimes V^{+}\,|\,h_{n}v=\delta_{n,0}\nu^{+}(\mu)(h)v
\text{ for any }h\in\mf{t}^{+}\text{ and }n\geq0\}.
\end{align*}
\item When $V=V_{\sf sc}$, we define
a $\mathbb{Z}/2\mathbb{Z}$-graded weak $V_{\sf af}$-module\footnote{
When $V=\mathcal{C}_{k}(\mf{g})$, we use Proposition \ref{strong_FST}.} by
\begin{align*}
\Omega^{-}(\mu)(M):=\{v\in M\otimes V^{-}\,|\,h_{n}v=\delta_{n,0}\nu^{-}(\mu)(h)v
\text{ for any }h\in\mf{t}^{-}\text{ and }n\geq0\}.
\end{align*}
\end{enumerate}
For a morphism $f\colon M^{1}\to M^{2}$ of $\mathscr{C}_{[\lambda]}(V)$, we define
a linear map
$$\Omega^{\pm}(\mu)(f)\colon \Omega^{\pm}(\mu)(M^{1})\to\Omega^{\pm}(\mu)(M^{2})$$ by
the restriction of $f\otimes{\sf id}_{V^{\pm}}$ to the subspace $\Omega^{\pm}(\mu)(M^{1})$
of $M^{1}\otimes V^{\pm}$.
\end{dfn}

We denote the purely even 
Heisenberg Fock $\mathcal{H}^{\pm}$-module
of highest weight $\mu^{\pm}\in(\mf{t}^{\pm})^{*}$
by $\mathcal{H}^{\pm}_{\mu^{\pm}}$
and its even highest weight vector by $\ket{\mu^{\pm}}$.

\begin{lem}\label{decomposition}
For $V=V_{\sf af}$, there exists a unique isomorphism
\begin{equation}\label{branch1}
\bigoplus_{\mu\in[\lambda]}\Omega^{+}(\mu)(M)
\otimes\mathcal{H}^{+}_{\nu^{+}(\mu)}
\xrightarrow{\simeq}\iota_{+}^{*}(M\otimes V^{+})
\end{equation}
of $\mathbb{Z}/2\mathbb{Z}$-graded weak $V_{\sf sc}\otimes\mathcal{H}^{+}$-modules
such that $v\otimes\ket{\nu^{+}(\mu)}\mapsto v$ for any $v\in\Omega^{+}(\mu)(M)$.
Similarly, for $V=V_{\sf sc}$, 
there exists a unique isomorphism
\begin{equation}\label{branch2}
\bigoplus_{\mu\in[\lambda]}\Omega^{-}(\mu)(M)\otimes\mathcal{H}^{-}_{\nu^{-}(\mu)}
\xrightarrow{\simeq}\iota_{-}^{*}(M\otimes V^{-})
\end{equation}
of $\mathbb{Z}/2\mathbb{Z}$-graded weak $V_{\sf af}\otimes\mathcal{H}^{-}$-modules
such that $v\otimes\ket{\nu^{-}(\mu)}\mapsto v$ for any $v\in\Omega^{-}(\mu)(M)$.
\end{lem}

\begin{proof}
Since $\iota_{\pm}^{*}(M\otimes V^{\pm})$ is completely reducible 
as a $\mathbb{Z}/2\mathbb{Z}$-graded weak
$\mathcal{H}^{\pm}$-module by the same argument as in \cite[\S1.7]{frenkel1989vertex}
and $\Omega^{\pm}(\mu)(M)=\{0\}$ holds for any $\mu\in\mf{t}_{V}^{*}\setminus[\lambda]$,
we get the above isomorphisms.
\end{proof}

\begin{prp}\label{Funct}
We have the following.
\begin{enumerate}
\item
Let $\lambda\in\mf{h}^{*}$ and $\mu\in[\lambda]$. 
Then the assignment $\Omega^{+}(\mu)$ 
gives rise to a $\mathbb{C}$-linear functor
\begin{equation*}
\Omega^{+}(\mu)\colon
\mathscr{C}_{[\lambda]}(V_{\sf af})
\to\mathscr{C}_{[\mu_{\sf sc}]}(V_{\sf sc}),
\end{equation*}
where $\mu_{\sf sc}\in(\mf{h}_{\sf sc})^{*}$ is defined by
$$\mu_{\sf sc}(J_{\alpha}):=\mu\left(k^{-1}H_{\alpha}\right)$$
for $\alpha\in\Delta^{+}$.
\item
Let $\lambda\in(\mf{h}_{\sf sc})^{*}$
and $\mu\in[\lambda]$. 
Then the assignment $\Omega^{-}(\mu)$
gives rise to a $\mathbb{C}$-linear functor
$$\Omega^{-}(\mu)\colon
\mathscr{C}_{[\lambda]}(V_{\sf sc})
\to\mathscr{C}_{[\mu_{\sf af}]}(V_{\sf af}),
$$
where $\mu_{\sf af}\in\mf{h}^{*}$ is defined by
$$\mu_{\sf af}(H_{\alpha}):=
\mu\left(kJ_{\alpha}\right)$$
for any $\alpha\in\Pi$.
\end{enumerate}
\end{prp}

\begin{proof}
Since the proof of (1) and that of (2) are simliar, we only verify (1) here.
Let $M$ be an object of $\mathscr{C}_{[\lambda]}(V_{\sf af})$. 
As the functoriality is obvious, 
it suffices to show that $\Omega^{+}(\mu)(M)$ is an object of
$\mathscr{C}_{[\mu_{\sf sc}]}(V_{\sf sc})$.
By the definition of $\Omega^{+}(\mu)(M)$, the operator $H_{\alpha,0}^{+}$
acts on $\Omega^{+}(\mu)(M)$ as scalar 
$\nu^{+}(\mu)(H_{\alpha}^{+})=\mu(H_{\alpha})$.
Therefore $J_{\alpha,0}$ acts on
$\Omega^{+}(\mu)(M)$ as scalar 
$\mu_{\sf sc}(J_{\alpha})$ modulo $\mathbb{Z}$
and every $\mf{h}_{\sf sc}$-weight of
$\Omega^{+}(\mu)(M)$ lies in the coset $[\mu_{\sf sc}]$.
The other conditions of $\mathscr{C}_{[\mu_{\sf sc}]}(V_{\sf sc})$
are easily verified by the corresponding conditions
of $\mathscr{C}_{[\lambda]}(V_{\sf af})$.
This completes the proof.
\end{proof}

The next lemma is obvious by the definition.

\begin{lem}\label{sc_af}
For any $\lambda\in\mf{h}^{*}$, we have
$(\lambda_{\sf sc})_{\sf af}=\lambda$.
\end{lem}

\subsection{The sublattice $K$}\label{lattice_K}

We define a $\mathbb{Z}$-linear map 
$g_{\sf af}^{+}\colon L^{+}\to Q$ by
\begin{equation}\label{g_plus}
g_{\sf af}^{+}(\xi):=\sum_{\alpha\in\Delta^{+}}\langle\xi,\alpha^{+}\rangle\alpha
\end{equation}
for $\xi\in L^{+}$
and regard its kernel $K$ as a positive-definite sublattice
of $L^{+}$.
Since we have $g_{\sf af}^{+}\circ f^{+}_{\sf af}={\sf id}_{Q}$,
it is easy to verify that
the set
$$\big\{\xi(\alpha):=\alpha^{+}-f^{+}_{\sf af}(\alpha)\in L^{+}\,
\big|\,\alpha\in\Delta^{+}\setminus\Pi\big\}$$
forms a $\mathbb{Z}$-basis of $K$.

\begin{lem}
Let $V_{K}$ be the lattice vertex superalgebra  
associated with the lattice $K$. 
We identify $V_{K}$ with
the vertex subsuperalgebra of $V_{\sf af}\otimes V^{+}$
generated by
$\{\mathbf{1}_{\sf af}\otimes\mathrm{e}^{\xi}\,|\,\xi\in K\}$.
Then $V_{K}$ is contained in $V_{\sf sc}$.
In addition, we have
\begin{equation}\label{simple_reduction}
\widetilde{\xi}(\alpha):=\mathbf{1}_{\sf af}\otimes\xi(\alpha)_{-1}\mathbf{1}_{V^{+}}
=J_{\alpha}-\sum_{\beta\in\Pi}\alpha(\beta^{*})J_{\beta}\in\mf{h}_{\sf sc}
\end{equation}
for any $\alpha\in\Delta^{+}$.
\end{lem}

\begin{proof}
Since $\langle\alpha_{\sf f},\xi(\beta)\rangle=0$ holds for any $\alpha\in\Pi$
and $\beta\in\Delta^{+}\setminus\Pi$,
the set of generators is contained in $V_{\sf sc}$.
The equality (\ref{simple_reduction}) follows from
$$H_{\alpha}-\sum_{\beta\in\Pi}\alpha(\beta^{*})H_{\beta}=0.$$
\end{proof}

With the help of the vertex subsuperalgebra $V_{K}$ of $V_{\sf sc}$,
we obtain the `converse' of Lemma \ref{sc_af} as follows.

\begin{prp}\label{converse}
Let $\lambda\in(\mf{h}_{\sf sc})^{*}$.
Assume that
$\mathscr{C}_{[\lambda]}(V_{\sf sc})$
contains at least one non-zero object.
Then $[(\mu_{\sf af})_{\sf sc}]=[\lambda]$ holds for any $\mu\in[\lambda]$.
\end{prp}

\begin{proof}
Take $\mu\in[\lambda]$.
When $\mf{g}=\mf{sl}_{2}$, we have $(\mu_{\sf af})_{\sf sc}=\mu$
and $[(\mu_{\sf af})_{\sf sc}]=[\mu]=[\lambda]$.
Now we assume that $\mf{g}\neq\mf{sl}_{2}$.
It suffices to show that 
$\mu(J_{\alpha})
-(\mu_{\sf af})_{\sf sc}(J_{\alpha})$ lies in $\mathbb{Z}$
for any $\alpha\in\Delta^{+}$.
Since $K$ is non-trivial and $V_{K}$ is regular (see \cite[Theorem 3.16]{dong1997regularity}),
any non-zero object of
$\mathscr{C}_{[\lambda]}(V_{\sf sc})$
decomposes into a non-empty direct sum
of simple ordinary $V_{K}$-modules.
Therefore we have $\mu\big(\widetilde{\xi}(\alpha)\big)\in\mathbb{Z}$
for any $\alpha\in\Delta^{+}$.
On the other hand, by the definition of $(\mu_{\sf af})_{\sf sc}$ and (\ref{simple_reduction}),
we also have 
$(\mu_{\sf af})_{\sf sc}\big(\widetilde{\xi}(\alpha)\big)=0$.
Then, by using (\ref{simple_reduction}) 
and $\mu(J_{\beta})=(\mu_{\sf af})_{\sf sc}(J_{\beta})$ for any $\beta\in\Pi$, we obtain 
$\mu(J_{\alpha})-(\mu_{\sf af})_{\sf sc}(J_{\alpha})\in\mathbb{Z}$
for any $\alpha\in\Delta^{+}$.
We thus complete the proof.
\end{proof}

\subsection{Main Theorem}

Our main result is as follows:

\begin{thm}\label{MainThm}
Let $k\in\mathbb{C}\setminus\{0,-h^{\vee}\}$.
\begin{enumerate}
\item
For any $[\lambda]\in\mf{h}^{*}/Q$
and $\mu\in[\lambda]$, the following $\mathbb{C}$-linear functors
\begin{align*}
&\Omega^{+}(\mu)\colon
\mathscr{C}_{[\lambda]}(V_{\sf af})
\to\mathscr{C}_{[\mu_{\sf sc}]}(V_{\sf sc}),\\
&\Omega^{-}(\mu_{\sf sc})\colon
\mathscr{C}_{[\mu_{\sf sc}]}(V_{\sf sc})
\to\mathscr{C}_{[\lambda]}(V_{\sf af})
\end{align*}
are mutually quasi-inverse to each other.
\item
For any $\lambda\in\mf{h}^{*}$ and $\gamma\in Q$,
there exist an isomorphism
$${\sf U}^{\gamma}_{\sf sc}\colon\mathscr{C}_{[\lambda_{\sf sc}]}(V_{\sf sc})
\xrightarrow{\simeq}\mathscr{C}_{[(\lambda-\gamma)_{\sf sc}]}(V_{\sf sc})$$
of categories and an isomorphism
$$\Omega^{+}(\lambda-\gamma)\simeq{\sf U}^{\gamma}_{\sf sc}\circ
\Omega^{+}(\lambda)$$
of functors from $\mathscr{C}_{[\lambda]}(V_{\sf af})$
to $\mathscr{C}_{[(\lambda-\gamma)_{\sf sc}]}(V_{\sf sc})$.
\item
For any $\lambda\in(\mf{h}_{\sf sc})^{*}$ and $\gamma\in Q_{\sf sc}$,
there exist an isomorphism
$${\sf U}^{\gamma}_{\sf af}\colon\mathscr{C}_{[\lambda_{\sf af}]}(V_{\sf af})
\xrightarrow{\simeq}\mathscr{C}_{[(\lambda+\gamma)_{\sf af}]}(V_{\sf af})$$
of categories and an isomorphism
$$\Omega^{-}(\lambda+\gamma)\simeq{\sf U}^{\gamma}_{\sf af}\circ
\Omega^{-}(\lambda)$$
of functors from $\mathscr{C}_{[\lambda]}(V_{\sf sc})$ to 
$\mathscr{C}_{[(\lambda+\gamma)_{\sf af}]}(V_{\sf af})$.
\end{enumerate}
\end{thm}

\begin{rem}
Following the notation and terminology of \cite{brundan2017monoidal},
we give a comment on supercategory structures.
Since the category $V\text{-\sf gMod}$ 
is a $\mathbb{C}$-linear supercategory (see Lemma \ref{Supercategory}
for the definition),
so is the full subcategory $\mathscr{C}_{[\lambda]}(V)$.
It is clear by the definition that $\Omega^{+}(\mu)$ and $\Omega^{-}(\mu_{\sf sc})$
are \emph{superfunctors} in the sense of \cite[Definition 1.1]{brundan2017monoidal}.
In addition, by the proof in the next section (see Remark \ref{superequiv1} and \ref{superequiv2}), 
the underlying categories
$\underline{\mathscr{C}_{[\lambda]}(V_{\sf af})}$
and
$\underline{\mathscr{C}_{[\mu_{\sf sc}]}(V_{\sf sc})}$
in the sense of \cite[Definition 1.1]{brundan2017monoidal}
are also categorical equivalent to each other.
\end{rem}

\section{Proof of Main Theorem}\label{ProofSection}

In this subsection, we fix $\lambda\in\mf{h}^{*}$ and $\mu\in[\lambda]$. 

\subsection{Faithfulness of $\Omega^{+}(\mu)$}

In this subsection, without using Proposition \ref{strong_FST},
we prove that 
$$\Omega^{-}(\mu_{\sf sc})\circ\Omega^{+}(\mu)\colon
\mathscr{C}_{[\lambda]}\big(V^{k}(\mf{g})\big)\to\mathscr{C}_{[\lambda]}\big(V^{k}(\mf{g})\big)
$$ is 
naturally isomorphic to the identity functor.
As corollaries, we obtain 
the proof of Proposition \ref{strong_FST}
and the faithfulness of $\Omega^{+}(\mu)$.

\subsubsection{Character formula}

Similarly to (\ref{g_plus}), 
we define $\mathbb{Z}$-linear maps
$g_{\sf af}^{-}\colon L^{-}\to Q$ 
and $g_{\sf sc}^{\pm}\colon L^{\pm}\to Q_{\mathcal{C}^{k}(\mf{g})}$ by
\begin{align*}
&g_{\sf af}^{-}(\zeta):=\sum_{\alpha\in\Pi}\langle\zeta,\alpha^{-}\rangle\alpha,\\
&g_{\sf sc}^{+}(\xi)(J_{\beta}^{*}):=\langle\xi,\beta^{+}\rangle\text{ for }\beta\in\Delta^{+},\\
&g_{\sf sc}^{-}(\zeta)(J_{\beta}^{*}):=
\begin{cases}
\langle\zeta,\beta^{-}\rangle &\text{ if }\beta\in\Pi,\\
0 & \text{ if }\beta\in\Delta^{+}\setminus\Pi
\end{cases}
\end{align*}
for $\xi\in L^{+}$ and $\zeta\in L^{-}$.

\begin{prp}\label{ChFormula}
Let $M$ be an object of $\mathscr{C}_{[\lambda]}\big(V^{k}(\mf{g})\big)$.
Then we have
\begin{equation}\label{CF1}
\ch\big(\Omega^{+}(\mu)(M)\big)
=\sum_{\xi\in L^{+}}
\frac{s^{\mu+g_{\sf af}^{+}(\xi)}_{M}(q)
q^{\frac{1}{2}\langle\xi,\xi\rangle-\Delta_{\mu}^{+}}}
{\eta(q)^{N-\ell}}\mathrm{e}^{\mu_{\sf sc}+g_{\sf sc}^{+}(\xi)}
\end{equation}
and
\begin{equation}\label{CF2}
\ch\big(\Omega^{-}(\mu_{\sf sc})\circ\Omega^{+}(\mu)(M)\big)
=\ch(M),
\end{equation}
where
$\Delta_{\mu}^{+}:=\frac{1}{k+h^{\vee}}\frac{(\mu,\mu)}{2}$
is the lowest $L_{0}^{\sf b}$-eigenvalue of $\mathcal{H}^{+}_{\nu^{+}(\mu)}$.
\end{prp}

\begin{proof}
We first prove (\ref{CF1}).
One can verify that the restriction of the isomorphism (\ref{branch1})
for $V=V^{k}(\mf{g})$ gives
an even linear isomorphism
$$\Omega^{+}(\mu)(M)\otimes\mathcal{H}^{+}_{\nu^{+}(\mu)}\simeq
\bigoplus_{\xi\in L^{+}}M_{\mu+g_{\sf af}^{+}(\xi)}\otimes
V^{+}_{\xi},$$
where
$V^{+}_{\xi}:=\big\{v\in V^{+}\,\big|\,\xi'_{0}v=\langle\xi,\xi'\rangle v\text{ for any }
\xi'\in L^{+}\big\}.$
It follows by direct calculation that
$M_{\mu+g_{\sf af}^{+}(\xi)}\otimes
V^{+}_{\xi}$ is the $\mf{t}_{\mathcal{C}^{k}(\mf{g})}$-eigenspace of eigenvalue 
$\mu_{\sf sc}+g_{\sf sc}^{+}(\xi)$,
the $(L^{\sf af}_{0}\otimes\id+\id\otimes L^{\sf f}_{0})$-graded dimension of 
$M_{\mu+g_{\sf af}^{+}(\xi)}\otimes
V^{+}_{\xi}$ is given by
$s^{\mu+g_{\sf af}^{+}(\xi)}_{M}(q)q^{\frac{1}{2}\langle\xi,\xi\rangle}\eta(q)^{-N}$,
and the $L^{\sf b}_{0}$-graded dimension of 
$\mathcal{H}^{+}_{\nu^{+}(\mu)}$ is given by $q^{\Delta_{\mu}^{+}}\eta(q)^{-\ell}$.
By combining them, we obtain the required formula (\ref{CF1}).

Next we verify (\ref{CF2}). 
For an arbitrary object $\mathcal{M}$ in 
$\mathscr{C}_{[\mu_{\sf sc}]}\big(\mathcal{C}^{k}(\mf{g})\big)$,
in the same way as above,
we can verify that the restriction of the isomorphism (\ref{branch2}) 
for $V=\mathcal{C}^{k}(\mf{g})$ gives
$$\Omega^{-}(\mu_{\sf sc})(\mathcal{M})\otimes
\mathcal{H}_{\nu^{-}(\mu_{\sf sc})}^{-}\simeq\bigoplus_{\zeta\in L^{-}}
\mathcal{M}_{\mu_{\sf sc}-g_{\sf sc}^{-}(\zeta)}\otimes V^{-}_{\zeta}$$
and
\begin{equation}\label{CF3}
\ch\big(\Omega^{-}(\mu_{\sf sc})(\mathcal{M})\big)
=\sum_{\zeta\in L^{-}}
\frac{s^{\mu_{\sf sc}-g_{\sf sc}^{-}(\zeta)}_{\mathcal{M}}(q)
q^{\frac{1}{2}\langle\zeta,\zeta\rangle+\Delta_{\mu}^{+}}}
{\eta(q)^{-N+\ell}}\mathrm{e}^{\mu-g_{\sf af}^{-}(\zeta)},
\end{equation}
where
$V^{-}_{\zeta}:=\big\{v\in V^{-}\,\big|\,\zeta'_{0}v=\langle\zeta,\zeta'\rangle v\text{ for any }
\zeta'\in L^{-}\big\}.$
On the other hand, by some computation, 
we obtain
\begin{equation}\label{CFlemma}
\sum_{(\xi,\zeta)\in S(\gamma)}
q^{\frac{1}{2}\langle\xi,\xi\rangle+\frac{1}{2}\langle\zeta,\zeta\rangle}
s^{\mu+g_{\sf af}^{+}(\xi)}_{M}(q)
=s^{\mu+\gamma}_{M}(q)
\end{equation}
for any $\gamma\in Q$, where 
$$S(\gamma):=\big\{(\xi,\zeta)\in L^{+}\times L^{-}\,\big|\,
\gamma=-g_{\sf af}^{-}(\zeta),\ g_{\sf sc}^{+}(\xi)=-g_{\sf sc}^{-}(\zeta)\big\}.$$
Then, by using (\ref{CF3}) for $\mathcal{M}=\Omega^{+}(\mu)(M)$ 
together with (\ref{CF1}) and (\ref{CFlemma}),
we can derive the character formula (\ref{CF2}).
\end{proof}

\subsubsection{Twisted embedding}

Let $M$ be an object of $\mathscr{C}_{[\lambda]}\big(V^{k}(\mf{g})\big)$.
We define an even linear operator $\mathcal{H}$ on $M\otimes V^{+}\otimes V^{-}$ by
$$\mathcal{H}:=\sum_{\alpha\in\Pi}\sum_{n=1}^{\infty}
\frac{1}{n}\alpha^{*}_{n}\otimes
\big(\alpha_{-n}^{+}\otimes\id
+\id\otimes\alpha_{-n}^{-}\big).$$
Then, by the condition (2) of Definition \ref{harishchandra},
the formal sum ${\sf exp}(\mathcal{H})$ defines 
an even linear automorphism on
$M\otimes V^{+}\otimes V^{-}$.

\begin{lem}\label{tw_emb}
The assignment
\begin{equation}\label{assign}
v\otimes\ket{\nu^{+}(\mu)}\otimes
\ket{\nu^{-}(\mu_{\sf sc})}\mapsto
\widetilde{v}:={\sf exp}(\mathcal{H})
\big(v\otimes\mathrm{e}^{f^{+}_{\sf af}(\gamma)}
\otimes\mathrm{e}^{f^{-}_{\sf af}(\gamma)}\big)
\end{equation}
for $v\in M_{\mu+\gamma}$ and $\gamma\in Q$
uniquely extends to an injective morphism
\begin{equation}\label{twisted_emb}
\iota_{M}^{\mu}\colon M\otimes
\mathcal{H}^{+}_{\nu^{+}(\mu)}
\otimes\mathcal{H}^{-}_{\nu^{-}(\mu_{\sf sc})}\to 
\iota_{\sf af}^{*}(M\otimes V^{+}\otimes V^{-})
\end{equation}
of $\mathbb{Z}/2\mathbb{Z}$-graded weak $V^{k}(\mf{g})\otimes
\mathcal{H}^{+}\otimes\mathcal{H}^{-}$-modules.
\end{lem}

\begin{proof}
First we prove that the assignment (\ref{assign})
defines a $\mathbb{Z}/2\mathbb{Z}$-graded weak $V^{k}(\mf{g})$-module homomorphism
from $M\otimes\ket{\nu^{+}(\mu)}\otimes
\ket{\nu^{-}(\mu_{\sf sc})}$
to $\iota_{\sf af}^{*}(M\otimes V^{+}\otimes V^{-}).$
It suffices to prove that
\begin{align}\label{twist}
\widetilde{X}_{\alpha,m}\widetilde{v}
={\sf exp}(\mathcal{H})\big(X_{\alpha,m}v\otimes\mathrm{e}^{f^{+}_{\sf af}(\gamma+\alpha)}
\otimes\mathrm{e}^{f^{-}_{\sf af}(\gamma+\alpha)}\big)
\end{align}
for any $\alpha\in\Pi\sqcup(-\Pi)$, $X_{\alpha}\in\mf{g}_{\alpha}$, and $m\in\mathbb{Z}$.
By direct calculation, we obtain
\begin{align*}
&\ad(\alpha^{*}_{n}\otimes \alpha_{-n}^{+})^{N}
(X_{\beta,p}\otimes e^{\beta^{+}}_{q})\\
&=\delta_{\alpha,\beta}\sum_{N'=0}^{N}
\binom{N}{N'}
\big(X_{\alpha,p+(N-N')n}\otimes(\alpha_{-n}^{+})^{N-N'}\big)
\circ\big((\alpha_{n}^{*})^{N'}\otimes e^{\alpha^{+}}_{q-N'n}\big)\\
&\ad(\alpha^{*}_{n}\otimes\alpha_{-n}^{-})^{N}
(X_{\beta,p}\otimes e^{\beta^{-}}_{r})\\
&=\delta_{\alpha,\beta}\sum_{N'=0}^{N}
(-1)^{N'}\binom{N}{N'}
\big(X_{\alpha,p+(N-N')n}\otimes (\alpha_{-n}^{-})^{N-N'}\big)
\circ\big((\alpha^{*})^{N'}\otimes e^{\alpha^{-}}_{r-N'n}\big)
\end{align*}
for any $\alpha,\beta\in\Pi$ and $p,q,r\in\mathbb{Z}$.
By using them, for $\alpha\in\Pi$, we have
\begin{align*}
\Phi&:={\sf exp}(-\mathcal{H})\,
\widetilde{X}_{\alpha,m}\,{\sf exp}(\mathcal{H})
\big(v\otimes\mathrm{e}^{f^{+}_{\sf af}(\gamma)}
\otimes\mathrm{e}^{f^{-}_{\sf af}(\gamma)}\big)\\
&=\sum_{p+q+r=m-2}\sum_{N=0}^{\infty}\frac{\ad(-\mathcal{H})^{N}}{N!}
(X_{\alpha,p}\otimes e^{\alpha^{+}}_{q}
\otimes e^{\alpha^{-}}_{r})
\big(v\otimes\mathrm{e}^{f^{+}_{\sf af}(\gamma)}
\otimes\mathrm{e}^{f^{-}_{\sf af}(\gamma)}\big)\\
&=X_{\alpha,m}v\otimes\mathrm{e}^{f^{+}_{\sf af}(\gamma+\alpha)}
\otimes\mathrm{e}^{f^{-}_{\sf af}(\gamma+\alpha)}\\
&+
\sum_{p+q+r=m-2}\sum_{N=1}^{\infty}\frac{(-1)^{N}}{N!}
\sum_{N'=0}^{N}
\binom{N}{N'}
\sum_{n=1}^{\infty}
X_{\alpha,p+(N-N')n}(\alpha^{*}_{n})^{N'}v\otimes\Psi_{N,N',n;q,r}
\end{align*}
where
\begin{align*}
\Psi_{N,N',n;q,r}
&:=(\alpha^{+}_{-n})^{N-N'}e^{\alpha^{+}}_{q-N'n}
\mathrm{e}^{f^{+}_{\sf af}(\gamma)}\otimes
e^{\alpha^{-}}_{r}\mathrm{e}^{f^{-}_{\sf af}(\gamma)}\\
&\hspace{2cm}+(-1)^{N'}e^{\alpha^{+}}_{q}
\mathrm{e}^{f^{+}_{\sf af}(\gamma)}\otimes
(\alpha^{-}_{-n})^{N-N'}e^{\alpha^{-}}_{r-N'n}\mathrm{e}^{f^{-}_{\sf af}(\gamma)}.
\end{align*}
It is clear that $\Phi$ lies in $V^{k}(\mf{g})\otimes 
V^{+}_{f^{+}_{\sf af}(\gamma+\alpha)}\otimes 
V^{-}_{f^{-}_{\sf af}(\gamma+\alpha)}$.
By straightforward computations, we have
$(\id\otimes(\alpha')^{+}_{n}\otimes\id)\Phi=
(\id\otimes\id\otimes(\alpha'')^{-}_{n})\Phi=0$
for any $\alpha'\in\Delta^{+}$, $\alpha''\in\Pi$, and $n\in\mathbb{Z}_{>0}$.
Therefore, by the uniquness of singular vector in the Heisenberg Fock module
$V^{+}_{f^{+}_{\sf af}(\gamma+\alpha)}\otimes 
V^{-}_{f^{-}_{\sf af}(\gamma+\alpha)}$,
there exists $\Phi_{\sf af}\in V^{k}(\mf{g})$ such that
$\Phi=\Phi_{\sf af}\otimes e^{f^{+}_{\sf af}(\gamma+\alpha)}\otimes 
e^{f^{-}_{\sf af}(\gamma+\alpha)}.$
Then, by the explicit form of $\Psi$, 
we conclude that
$\Phi_{\sf af}=X_{\alpha,m}v$.
This proves the formula (\ref{twist}) for $\alpha\in\Pi$. Since the proof 
for $\alpha\in-\Pi$ is the same,
we omit it.

Next we prove that the assignment (\ref{assign})
uniquely extends to (\ref{twisted_emb}).
By a direct computation, we obtain
\begin{align*}
H_{\alpha,n}^{+}\widetilde{v}
=\begin{cases}
0 & \text{ if }n>0,\\
\nu^{+}(\mu)(H_{\alpha}^{+})\widetilde{v} & \text{ if }n=0,
\end{cases}
\ \ H_{\beta,n}^{-}\widetilde{v}
=\begin{cases}
0 & \text{ if }n>0,\\
\nu^{-}(\mu_{\sf sc})(H_{\beta}^{-})\widetilde{v} & \text{ if }n=0
\end{cases}
\end{align*}
for $\alpha\in\Pi$ and $\beta\in\Delta^{+}$.
Hence (\ref{assign}) uniquely extends to (\ref{twisted_emb}). 

At last, as the injectivity of (\ref{assign}) follows from 
the bijectivity of the operator ${\sf exp}(\mathcal{H})$
on $M\otimes V^{+}\otimes V^{-}$,
the induced homomorphism (\ref{twisted_emb}) is also injective.
\end{proof}

\subsubsection{Natural isomorphism}

\begin{prp}\label{natural}
Let $M$ be an object of $\mathscr{C}_{[\lambda]}\big(V^{k}(\mf{g})\big)$.
Then the injective homomorphism (\ref{twisted_emb}) gives rise to
an isomorphism
\begin{equation}\label{naturalisom}
\mathscr{F}^{\mu}_{M}\colon M\xrightarrow{\simeq}
\Omega^{-}(\mu_{\sf sc})\circ\Omega^{+}(\mu)(M);\ 
v\mapsto \iota^{\mu}_{M}\big(v\otimes\ket{\nu^{+}(\mu)}\otimes\ket{\nu^{-}(\mu_{\sf sc})}\!\big)
\end{equation}
of $\mathbb{Z}/2\mathbb{Z}$-graded weak $V^{k}(\mf{g})$-modules,
which is natural in $M$.
\end{prp}

\begin{proof}
First we prove that $\mathscr{F}^{\mu}_{M}$ is a $\mathbb{Z}/2\mathbb{Z}$-graded 
weak $V^{k}(\mf{g})$-module isomorphism.
Since $\mathscr{F}^{\mu}_{M}$ is a well-defined injective
morphism of $\mathbb{Z}/2\mathbb{Z}$-graded 
weak $V^{k}(\mf{g})$-modules by Lemma \ref{tw_emb}, it suffices to show that
$\mathscr{F}^{\mu}_{M}$ is surjective.
The surjectivity of $\mathscr{F}^{\mu}_{M}$
follows from its injectivity and the character formula (\ref{CF2}).

Next we prove that $\mathscr{F}^{\mu}_{M}$ is natural in $M$.
Let $f\colon M^{1}\to M^{2}$ be a morphism of $\mathscr{C}_{[\lambda]}\big(V^{k}(\mf{g})\big)$.
By using the explicit form (\ref{assign}) of $\mathscr{F}^{\mu}_{M}$, we have
\begin{align*}
\big(\Omega^{-}(\mu_{\sf sc})\circ\Omega^{+}(\mu)\big)(f)
\big(\mathscr{F}^{\mu}_{M^{1}}(v)\big)
&=(f\otimes\id_{V^{+}}\otimes\id_{V^{-}})
\big(\mathscr{F}^{\mu}_{M^{1}}(v)\big)\\
&=\mathscr{F}^{\mu}_{M^{2}}\big(f(v)\big)
\end{align*}
for any $v\in M^{1}$.
We thus complete the proof.
\end{proof}

\begin{rem}\label{superequiv1}
It is clear by the proof that the above natural isomorphism is an 
\emph{even supernatural transformation} in the sense of \cite[Definition 1.1]{brundan2017monoidal}.
\end{rem}

Now we give the proof of Proposition \ref{strong_FST}.

\begin{cor}\label{FSTisom}
Proposition \ref{strong_FST} holds.
\end{cor}

\begin{proof}
Since $f_{\sf FST}=\mathscr{F}^{0}_{V^{k}(\mf{g})}$,
it follows from Proposition \ref{natural}.
\end{proof}

By using Proposition \ref{strong_FST}, we obtain the following.

\begin{prp}\label{faithful}
The composed functor
$\Omega^{-}(\mu_{\sf sc})\circ\Omega^{+}(\mu)\colon
\mathscr{C}_{[\lambda]}(V_{\sf af})\to\mathscr{C}_{[\lambda]}(V_{\sf af})$
is naturally isomorphic to the identity functor of $\mathscr{C}_{[\lambda]}(V_{\sf af})$.
In particular, the functor $\Omega^{+}(\mu)\colon\mathscr{C}_{[\lambda]}(V_{\sf af})\to
\mathscr{C}_{[\mu_{\sf sc}]}(V_{\sf sc})$ is faithful.
\end{prp}

\begin{proof}
By Proposition \ref{strong_FST}, we can verify that Proposition \ref{natural}
remains true if we replace $V^{k}(\mf{g})$ by $L_{k}(\mf{g})$.
Therefore (\ref{naturalisom}) gives a desired natural isomorphism.
\end{proof}

\begin{rem}\label{difficulty}
When $\mf{g}=\mf{sl}_{2}$,
one can describe an explicit natural isomorphism
from the identity functor of $\mathscr{C}_{\mu_{\sf sc}}(V_{\sf sc})$
to $\Omega^{+}(\mu)\circ\Omega^{-}(\mu_{\sf sc})$
in the same way as \cite[\S5.3]{sato2016equivalences}.
On the other hand, when $\mf{g}\neq\mf{sl}_{2}$,
there exist certain ``inner automorphisms'' of $V_{\sf sc}$ induced by the lattice $K$
(see Lemma \ref{inner_symm})
and the same argument as \cite[\S5.3]{sato2016equivalences}
is not applicable.
\end{rem}

\subsection{Essentially surjectivity and fullness of $\Omega^{+}(\mu)$}\label{EssFull}

In this subsection, we prove that
$\Omega^{+}(\mu)\colon\mathscr{C}_{[\lambda]}(V_{\sf af})\to
\mathscr{C}_{[\mu_{\sf sc}]}(V_{\sf sc})$ is essentially surjective and full.

We first recall the theory of Li's $\Delta$-operators developed in \cite{li1997physics}.
Let $(\mathcal{V}=\mathcal{V}^{\bar{0}}\oplus\mathcal{V}^{\bar{1}},
Y,\mathbf{1},\omega)$ be a conformal vertex superalgebra.
We use the following notation introduced in \cite[(3.4)]{li1997physics}:
$$E^{\pm}(h,z):=
{\sf exp}\left(\sum_{j=1}^{\infty}\frac{h_{\pm j}}{j}z^{\mp j}\right)\in
{\sf End}(\mathcal{V})[\![z^{\mp1}]\!]$$
for $h\in\mathcal{V}^{\bar{0}}$.

\begin{lem}[\cite{li1997physics}]\label{LiDelta}
Let $(\mathcal{M},Y_{\mathcal{M}})$ be a 
$\mathbb{Z}/2\mathbb{Z}$-graded weak $\mathcal{V}$-module.
We assume that there exist 
$h\in\mathcal{V}^{\bar{0}}$ and $c\in\mathbb{C}\setminus\{0\}$ satisfying 
the following conditions:
\begin{enumerate}
\item $L_{n}h=\delta_{n,0}h$ holds for any $n\geq0$,
\item $h_{n}h=\delta_{n,1}c\mathbf{1}$ holds for any $n\geq0$,
\item the operator $h_{0}$ on $\mathcal{V}$ is diagonalizable and 
has only integral eigenvalues,
\item 
$E^{+}(-h,-z)v\in\mathcal{V}[z,z^{-1}]$ holds for any $v\in\mathcal{V}$.
\end{enumerate}
Then Li's $\Delta$-operator
$$\Delta(h,z):=z^{h_{0}}E^{+}(-h,-z)\in{\sf End}(\mathcal{V})[\![z,z^{-1}]\!]$$
lies in the group 
$G^{0}(\mathcal{V})$ introduced in
\cite[\S2]{li1997physics}.
In particular,
we can \emph{twist} the $\mathbb{Z}/2\mathbb{Z}$-graded
weak $\mathcal{V}$-module structure $Y_{\mathcal{M}}$ on $\mathcal{M}$
to obtain a new $\mathbb{Z}/2\mathbb{Z}$-graded weak $\mathcal{V}$-module 
$\mathcal{M}^{h}:=(\mathcal{M},Y_{\mathcal{M}^{h}})$ defined by
$$Y_{\mathcal{M}^{h}}=Y_{\mathcal{M}^{h}}(?,z):=Y_{\mathcal{M}}\big(\Delta(h,z)?;z\big)\colon
\mathcal{V}\to{\sf End}(\mathcal{M})[\![z,z^{-1}]\!].$$
\end{lem}

\begin{proof}
See \cite[Proposition 2.1 and 3.2]{li1997physics} for the proof.
\end{proof}

We now focus on the case of $\mathcal{V}=\mathcal{C}^{k}(\mf{g})\otimes\mathcal{H}^{+}$.

\begin{lem}\label{K_twist}
Let $\mathcal{M}$ be an object of $\mathscr{C}_{[\mu_{\sf sc}]}\big(\mathcal{C}^{k}(\mf{g})\big)$.
We consider the direct sum
\begin{equation}\label{decompo}
\widetilde{\mathcal{M}}=\bigoplus_{\gamma\in Q}
\widetilde{\mathcal{M}}(\gamma)
:=\bigoplus_{\gamma\in Q}
\boldsymbol{\Pi}^{\langle f^{+}_{\sf af}(\gamma),f^{+}_{\sf af}(\gamma)\rangle}
\Big(\mathcal{M}\otimes\mathcal{H}_{\nu^{+}(\mu)}^{+}\Big)^{h(\gamma)}
\end{equation}
of $\mathbb{Z}/2\mathbb{Z}$-graded weak $\mathcal{C}^{k}(\mf{g})\otimes\mathcal{H}^{+}$-modules,
where 
$\boldsymbol{\Pi}$ is the $\mathbb{Z}/2\mathbb{Z}$-parity reversing functor and
$$h(\gamma):=\sum_{\alpha\in\Pi}\gamma(\alpha^{*})
\left(J_{\alpha}^{*}\otimes\mathbf{1}^{+}
-\mathbf{1}_{\sf sc}\otimes\frac{1}{k+h^{\vee}}H_{\alpha}^{+}\right)
\in\mathcal{C}^{k}(\mf{g})\otimes\mathcal{H}^{+}.$$
Then the $\mathbb{Z}/2\mathbb{Z}$-graded weak 
$\mathcal{C}^{k}(\mf{g})\otimes\mathcal{H}^{+}$-module structure
on $\widetilde{\mathcal{M}}$
can be extended to a $\mathbb{Z}/2\mathbb{Z}$-graded 
weak $V^{k}(\mf{g})\otimes V^{+}$-module structure 
on $\widetilde{\mathcal{M}}$.
\end{lem}

\begin{proof}
First we extend the $\mathbb{Z}/2\mathbb{Z}$-graded action of $V_{K}$ 
on $\widetilde{\mathcal{M}}$ to that of $V^{+}$.
Since we have
$\iota_{+}\big(h(\gamma)\big)
=\mathbf{1}_{\sf af}\otimes f^{+}_{\sf af}(\gamma)$,
the action of the Heisenberg vertex algebra 
generated by $\big\{h(\alpha)\,\big|\,\alpha\in\Pi\big\}$
on $\widetilde{\mathcal{M}}$
can be extended to that of $V_{f^{+}_{\sf af}(Q)}$
along with the $2$-cocycle $\epsilon^{+}$ 
(cf.\,\cite[Theorem 6.5.18]{lepowsky2012introduction}).
More precisely, we define 
a family of $\mathbb{Z}/2\mathbb{Z}$-homogeneous linear automorphisms 
$\{\mathbf{S}^{h(\gamma)}\,|\,\gamma\in Q\}$ on $\widetilde{\mathcal{M}}$ by
$$\mathbf{S}^{h(\gamma)}\Big|_{\widetilde{\mathcal{M}}(\gamma')}
:=\epsilon^{+}\big(f^{+}_{\sf af}(\gamma),f^{+}_{\sf af}(\gamma')\big)\id_{\mathcal{M}\otimes\mathcal{H}_{\nu^{+}(\mu)}^{+}}\colon\widetilde{\mathcal{M}}(\gamma')
\to\widetilde{\mathcal{M}}(\gamma'+\gamma)$$
and then the following $\mathbb{Z}/2\mathbb{Z}$-homogeneous mutually local fields\footnote{
The $\mathbb{Z}/2\mathbb{Z}$-parity of $\mathrm{e}^{f^{+}_{\sf af}(\gamma)}(z)$
is given by $\langle f^{+}_{\sf af}(\gamma),f^{+}_{\sf af}(\gamma)\rangle\ {\sf mod}\ 2\mathbb{Z}$.}
\begin{align}\label{VOexpression1}
\Big\{\mathrm{e}^{f^{+}_{\sf af}(\gamma)}(z):=\mathbf{S}^{h(\gamma)}
z^{h(\gamma)_{0}}
E^{-}\big(h(\gamma),z\big)
E^{+}\big(\!-h(\gamma),z\big)\,\Big|\,\gamma\in Q\Big\}
\end{align}
on $\widetilde{\mathcal{M}}$ generate the $\mathbb{Z}/2\mathbb{Z}$-graded
weak $V_{f^{+}_{\sf af}(Q)}$-module structure.
On the other hand, by direct calculation, we have
\begin{align*}
Y_{\widetilde{\mathcal{M}}}(\mathrm{e}^{\xi}\otimes\mathbf{1}^{+},z)
\Big|_{\widetilde{\mathcal{M}}(\gamma)}
=Y_{\mathcal{M}}(\mathrm{e}^{\xi},z)\otimes\id_{\mathcal{H}_{\nu^{+}(\mu)}^{+}}
\in{\sf End}\big(\widetilde{\mathcal{M}}(\gamma)\big)[\![z,z^{-1}]\!]
\end{align*}
for any $\gamma\in Q$.
By the general theory of regular lattice vertex superalgebra, 
there exists a family of $\mathbb{Z}/2\mathbb{Z}$-homogeneous
linear automorphisms $\{\mathbf{S}^{\xi}\,|\,\xi\in K\}$ of $\mathcal{M}$
such that
\begin{align}\label{VOexpression2}
Y_{\mathcal{M}}(\mathrm{e}^{\xi},z)
=\mathbf{S}^{\xi}z^{\xi_{0}}
E^{-}\left(\xi,z\right)
E^{+}\left(-\xi,z\right)
\end{align}
for any $\xi\in K$.
Then, by using the explicit description (\ref{VOexpression1})
and (\ref{VOexpression2}), one can verify that the following 
$\mathbb{Z}/2\mathbb{Z}$-homogeneous fields
$$\Big\{\mathrm{e}^{f^{+}_{\sf af}(\gamma)}(z),\ 
Y_{\widetilde{\mathcal{M}}}(\mathrm{e}^{\xi}\otimes\mathbf{1}^{+},z)\,\Big|\,
\gamma\in Q,\ \xi\in K\Big\}$$
are mutually local and generate the $\mathbb{Z}/2\mathbb{Z}$-graded
weak $V^{+}$-module structure on $\widetilde{\mathcal{M}}$
(cf.\,\cite[Theorem 3.2.10]{li1996localsys}).

Next, we further extend the above $V^{+}$-action on $\widetilde{\mathcal{M}}$ to
a $V^{k}(\mf{g})\otimes V^{+}$-action
which is compatible with the original 
$\mathcal{C}^{k}(\mf{g})\otimes\mathcal{H}^{+}$-action\footnote{
Since $\mathcal{C}^{k}(\mf{g})$ and $V^{+}\simeq\mathbf{1}_{\sf af}\otimes V^{+}$
generate the whole $V^{k}(\mf{g})\otimes V^{+}$,
such an extension is unique if it exists.}.
By using the residue product of mutually local fields 
(see e.g.\,\cite{kac97V}), 
we define the following mutually local even fields:
\begin{align*}
&X_{\alpha}(z):=
\epsilon^{+}\big(f^{+}_{\sf af}(\alpha),f^{+}_{\sf af}(-\alpha)\big)
Y_{\widetilde{\mathcal{M}}}(\Psi_{\alpha},z)
_{\langle f^{+}_{\sf af}(\alpha),f^{+}_{\sf af}(\alpha)\rangle-1}
\mathrm{e}^{-f^{+}_{\sf af}(\alpha)}(z),\\
&H_{\alpha}(z):=
Y_{\widetilde{\mathcal{M}}}\Big(\sum_{\beta\in\Delta^{+}}
(\alpha,\beta)J_{\beta}^{*}\otimes\mathbf{1}^{+}
+\mathbf{1}_{\sf sc}\otimes \frac{k}{k+h^{\vee}}H^{+}_{\alpha},z\Big)
\end{align*}
for $\alpha\in\Delta$,
where $\Psi_{\alpha}:=X_{\alpha,-1}\mathbf{1}_{\sf af}\otimes\mathrm{e}^{f^{+}_{\sf af}(\alpha)}
\otimes\mathbf{1}^{+}$.
Then, by using the Borcherds identity for mutually local fields,
we can verify that the set of mutually local fields
$\big\{X_{\alpha}(z),\,H_{\alpha}(z)\,
\big|\,\alpha\in\Delta\big\}$
generates a $\mathbb{Z}/2\mathbb{Z}$-graded weak
$V^{k}(\mf{g})$-module structure on $\widetilde{\mathcal{M}}$,
which commutes with the $V^{+}$-action on $\widetilde{\mathcal{M}}$.
\end{proof}

Since $L^{+}\simeq\mathbb{Z}^{N}$ is a positive-definite unimodular lattice, 
the canonical pairing
\begin{equation*}
{\sf Hom}_{V^{+}\text{-\sf Mod}}
(V^{+},\widetilde{\mathcal{M}})\otimes V^{+}\to\widetilde{\mathcal{M}};\ 
f\otimes v\mapsto f(v)
\end{equation*}
gives an isomorphism of 
$\mathbb{Z}/2\mathbb{Z}$-graded weak $V^{k}(\mf{g})\otimes V^{+}$-modules.

\begin{lem}\label{lift}
Let $\mathcal{M}$ and $\widetilde{\mathcal{M}}$ be as in Lemma \ref{K_twist}.
Then the $\mathbb{Z}/2\mathbb{Z}$-graded weak $V^{k}(\mf{g})$-module 
${\sf Hom}_{V^{+}\text{-\sf Mod}}(V^{+},\widetilde{\mathcal{M}})$
is an object of $\mathscr{C}_{[\lambda]}\big(V^{k}(\mf{g})\big)$.
\end{lem}

\begin{proof}
Since $\mathcal{M}$ is an object of $\mathscr{C}_{[\mu_{\sf sc}]}\big(\mathcal{C}^{k}(\mf{g})\big)$,
by the same computation as in Proposition \ref{ChFormula},
we can verify that
every $\mf{h}$-weight of ${\sf Hom}_{V^{+}\text{-\sf Mod}}
(V^{+},\widetilde{\mathcal{M}})$ lies in $[\lambda]$
and each string function of ${\sf Hom}_{V^{+}\text{-\sf Mod}}(V^{+},\widetilde{\mathcal{M}})$ 
has the lowest exponent in $q$.
Thus ${\sf Hom}_{V^{+}\text{-\sf Mod}}(V^{+},\widetilde{\mathcal{M}})$ 
is an object of $\mathscr{C}_{[\lambda]}\big(V^{k}(\mf{g})\big)$.
\end{proof}

\begin{prp}
The functor $\Omega^{+}(\mu)\colon\mathscr{C}_{[\lambda]}(V_{\sf af})\to
\mathscr{C}_{[\mu_{\sf sc}]}(V_{\sf sc})$ is essentially surjective and full.
\end{prp}

\begin{proof}
By Proposition \ref{faithful}, it suffices to prove the case of $(V_{\sf af}, V_{\sf sc})
=\big(V^{k}(\mf{g}),\mathcal{C}^{k}(\mf{g})\big)$.

We first prove that $\Omega^{+}(\mu)$
 is essentially surjective.
It suffices to show that $\mathcal{M}$ is isomorphic to
$\Omega^{+}(\mu)(M)$ for $M:={\sf Hom}_{V^{+}\text{-\sf Mod}}(V^{+},\widetilde{\mathcal{M}})$
as a $\mathbb{Z}/2\mathbb{Z}$-graded 
$\mathcal{C}^{k}(\mf{g})$-module.
By using Lemma \ref{LiDelta}, we obtain
a $\mathbb{Z}/2\mathbb{Z}$-graded $\mathcal{H}^{+}$-module isomorphism
\begin{equation*}
\Big(\mathcal{H}^{+}_{\nu^{+}(\mu)}\Big)^{\frac{1}{k+h^{\vee}}H_{\alpha}^{+}}
\xrightarrow{\simeq} \mathcal{H}^{+}_{\nu^{+}(\mu+\alpha)};\,
\ket{\nu^{+}(\mu)}\mapsto\ket{\nu^{+}(\mu+\alpha)}
\end{equation*}
for $\alpha\in\Pi$.
Hence we can rewrite the $\mathbb{Z}/2\mathbb{Z}$-graded
decomposition (\ref{decompo}) as
\begin{equation}\label{twisted_decomp}
\widetilde{\mathcal{M}}
\simeq\bigoplus_{\gamma\in Q}
\left(\boldsymbol{\Pi}^{\langle f^{+}_{\sf af}(\gamma),f^{+}_{\sf af}(\gamma)\rangle}\mathcal{M}^{h_{\sf sc}[\gamma]}
\right)\otimes\mathcal{H}_{\nu^{+}(\mu-\gamma)}^{+},
\end{equation}
where $h_{\sf sc}[\gamma]:=\sum_{\alpha\in\Pi}\gamma(\alpha^{*})J_{\alpha}^{*}$.
Therefore, by (\ref{branch1}) and Schur's lemma, 
we conclude that $\mathcal{M}\simeq\Omega^{+}(\mu)(M).$

Next we prove that $\Omega^{+}(\mu)$ is full.
Let $\mathcal{M}^{i}$ be an object of
$\mathscr{C}_{[\mu_{\sf sc}]}\big(\mathcal{C}^{k}(\mf{g})\big)$
and $\widetilde{\mathcal{M}}^{i}$ the extended $\mathbb{Z}/2\mathbb{Z}$-graded weak 
$V^{k}(\mf{g})\otimes V^{+}$-module as in Lemma \ref{K_twist} for $i\in\{1,2\}$.
Then an arbitrary morphism
$f\colon \mathcal{M}^{1}\to\mathcal{M}^{2}$ 
of weak $\mathcal{C}^{k}(\mf{g})$-modules extends to
a morphism of weak $V^{k}(\mf{g})\otimes V^{+}$-modules defined by
$$\widetilde{f}:=\bigoplus_{\gamma\in Q}f(\gamma)\colon
\widetilde{\mathcal{M}}^{1}\to\widetilde{\mathcal{M}}^{2} ,$$
where $f(\gamma)
:=f\otimes\id_{\mathcal{H}^{+}_{\nu^{+}(\mu)}}\colon
\widetilde{\mathcal{M}}^{1}(\gamma)\to\widetilde{\mathcal{M}}^{2}(\gamma)$.
Then, by the $\mathbb{Z}/2\mathbb{Z}$-graded linear isomorphisms
\begin{align*}
{\sf Hom}_{V^{k}(\mf{g})\otimes V^{+}\text{-\sf Mod}}
(\widetilde{\mathcal{M}}^{1},\widetilde{\mathcal{M}}^{2})
&\simeq{\sf Hom}_{V^{k}(\mf{g})\text{-\sf Mod}}(M^{1},M^{2})
\otimes{\sf Hom}_{V^{+}\text{-\sf Mod}}(V^{+},V^{+})\\
&\simeq{\sf Hom}_{V^{k}(\mf{g})\text{-\sf Mod}}(M^{1},M^{2}),
\end{align*}
we conclude that $f=\Omega^{+}(\mu)(\widetilde{f})$.
This completes the proof.
\end{proof}

\begin{rem}\label{superequiv2}
By the proof, the superfunctor $\Omega^{+}(\mu)$ is also
\emph{evenly dense} in the sense of \cite[Definition 1.1]{brundan2017monoidal}.
\end{rem}

\subsection{Spectral flow equivariance}\label{SpecFlowEquiv}
In this subsection, we prove (2) and (3) in Theorem \ref{MainThm}.

\subsubsection{Equivariance of $\Omega^{+}$}

\begin{lem}\label{U_sc}
Let $\lambda\in\mf{h}^{*}$ and $\gamma\in Q$.
For an object $\mathcal{M}$ and a morphism $f\colon\mathcal{M}^{1}\to\mathcal{M}^{2}$ 
of $\mathscr{C}_{[\lambda_{\sf sc}]}(V_{\sf sc})$,
we set
${\sf U}^{\gamma}_{\sf sc}(\mathcal{M})
:=\boldsymbol{\Pi}^{\langle f^{+}_{\sf af}(\gamma),f^{+}_{\sf af}(\gamma)\rangle}\mathcal{M}^{h_{\sf sc}[\gamma]}$
and
${\sf U}^{\gamma}_{\sf sc}(f):=f\colon{\sf U}^{\gamma}_{\sf sc}(\mathcal{M}^{1})
\to{\sf U}^{\gamma}_{\sf sc}(\mathcal{M}^{2}),$
respectively.
Then the assignment ${\sf U}^{\gamma}_{\sf sc}$ defines a $\mathbb{C}$-linear functor 
$${\sf U}^{\gamma}_{\sf sc}\colon\mathscr{C}_{[\lambda_{\sf sc}]}(V_{\sf sc})
\to\mathscr{C}_{[(\lambda-\gamma)_{\sf sc}]}(V_{\sf sc}),$$
which is a categorical isomorphism with the inverse 
${\sf U}^{-\gamma}_{\sf sc}$.
\end{lem}

\begin{proof}
By the two isomorphisms
(\ref{branch1}) and (\ref{twisted_decomp}) of
$\mathbb{Z}/2\mathbb{Z}$-graded $V_{\sf sc}\otimes\mathcal{H}^{+}$-modules,
the weak $V_{\sf sc}$-module
${\sf U}^{\gamma}_{\sf sc}(\mathcal{M})$ turns out to lie in 
$\mathscr{C}_{[(\lambda-\gamma)_{\sf sc}]}(V_{\sf sc})$.
The rest of the proof is straightforward and we omit it.
\end{proof}

\begin{rem}
It is known that ${\sf U}^{\gamma}_{\sf sc}$
induces the spectral flow automorphism
of the $\mathcal{N}=2$ superconformal algebra.
See \cite[Proposition 4.2]{hosono1991lie} for details.
\end{rem}

\begin{thm}\label{spf_equiv1}
For any $\lambda\in\mf{h}^{*}$ and $\gamma\in Q$,
we have an isomorphism
$$\Omega^{+}(\lambda-\gamma)\simeq{\sf U}^{\gamma}_{\sf sc}\circ\Omega^{+}(\lambda)$$
of functors from $\mathscr{C}_{[\lambda]}(V_{\sf af})$
to $\mathscr{C}_{[(\lambda-\gamma)_{\sf sc}]}(V_{\sf sc})$.
\end{thm}

\begin{proof}
By using (\ref{branch1}) and (\ref{twisted_decomp}), 
we can verify that
the natural isomorphism
$$M\otimes V^{+}\xrightarrow{\simeq} M\otimes
\boldsymbol{\Pi}^{\langle f^{+}_{\sf af}(\gamma),f^{+}_{\sf af}(\gamma)\rangle}(V^{+})^{f^{+}_{\sf af}(\gamma)};
m\otimes\mathbf{1}_{V^{+}}\mapsto
m\otimes\mathrm{e}^{-f^{+}_{\sf af}(\gamma)}$$
of $\mathbb{Z}/2\mathbb{Z}$-graded
weak $V_{\sf af}\otimes V^{+}$-modules
for an object $M$ of $\mathscr{C}_{[\lambda]}(V_{\sf af})$
gives rise to a natural isomorphism
$\Omega^{+}(\lambda-\gamma)(M)\xrightarrow{\simeq}
{\sf U}^{\gamma}_{\sf sc}\circ\Omega^{+}(\lambda)(M)$
of $\mathbb{Z}/2\mathbb{Z}$-graded weak $V_{\sf sc}$-modules.
\end{proof}

\subsubsection{Equivariance of $\Omega^{-}$}

\begin{lem}\label{inner_symm}
Let $\lambda\in(\mf{h}_{\sf sc})^{*}$ and $\xi\in K$.
For an object $\mathcal{M}$ 
and a morphism $f\colon\mathcal{M}^{1}\to\mathcal{M}^{2}$
of $\mathscr{C}_{[\lambda]}(V_{\sf sc})$, we set
${\sf U}^{\xi}_{\sf sc}(\mathcal{M}):=\boldsymbol{\Pi}^{\langle\xi,\xi\rangle}\mathcal{M}^{\xi}$
and ${\sf U}^{\xi}_{\sf sc}(f):=f\colon 
{\sf U}^{\xi}_{\sf sc}(\mathcal{M}^{1})\to{\sf U}^{\xi}_{\sf sc}(\mathcal{M}^{2})$,
respectively.
Then the assignment ${\sf U}^{\xi}_{\sf sc}$ defines a $\mathbb{C}$-linear functor
$${\sf U}^{\xi}_{\sf sc}\colon
\mathscr{C}_{[\lambda]}(V_{\sf sc})
\to\mathscr{C}_{[\lambda]}(V_{\sf sc}),$$
which is naturally isomorphic to 
the identity functor of $\mathscr{C}_{[\lambda]}(V_{\sf sc})$.
\end{lem}

\begin{proof}
By Proposition \ref{converse},
we may assume that $[\lambda]=[(\lambda_{\sf af})_{\sf sc}]$.
For an object $M$ of $\mathscr{C}_{[\lambda_{\sf af}]}(V_{\sf af})$,
in the same way as the proof of Theorem \ref{spf_equiv1},
the natural isomorphism
$$M\otimes V^{+}\xrightarrow{\simeq}
M\otimes \boldsymbol{\Pi}^{\langle\xi,\xi\rangle}(V^{+})^{\xi}; 
m\otimes\mathbf{1}_{V^{+}}\mapsto m\otimes\mathrm{e}^{-\xi}$$
of $\mathbb{Z}/2\mathbb{Z}$-graded weak $V_{\sf af}\otimes V^{+}$-modules
 gives rise to a natural isomorphism
$$\Omega^{+}(\lambda)(M)\xrightarrow{\simeq}
{\sf U}^{\xi}_{\sf sc}\circ\Omega^{+}(\lambda)(M)$$
of $\mathbb{Z}/2\mathbb{Z}$-graded weak $V_{\sf sc}$-modules.
Then, by Theorem \ref{MainThm}, we get the conclusion.
\end{proof}

\begin{thm}\label{spf_equiv2}
Let $\lambda\in(\mf{h}_{\sf sc})^{*}$ and $\gamma\in Q_{\sf sc}$. 
\begin{enumerate}
\item
For an object $M$ and a morphism $f\colon M^{1}\to M^{2}$ of
$\mathscr{C}_{[\lambda_{\sf af}]}(V_{\sf af})$, we set
${\sf U}^{\gamma}_{\sf af}(M):=\boldsymbol{\Pi}
^{\langle f^{+}_{\sf sc}(\gamma),f^{+}_{\sf sc}(\gamma)\rangle}M^{h_{\sf af}[\gamma]}$
and ${\sf U}^{\gamma}_{\sf af}(f):=f\colon {\sf U}^{\gamma}_{\sf af}(M^{1})\to
{\sf U}^{\gamma}_{\sf af}(M^{2})$, respectively,
where
$$h_{\sf af}[\gamma]:=\sum_{\alpha\in\Pi}\gamma(J_{\alpha}^{*})\alpha^{*}_{-1}\mathbf{1}_{\sf af},
\ \ f^{+}_{\sf sc}(\gamma):=\sum_{\beta\in\Delta^{+}}\gamma(J_{\beta}^{*})\beta^{+}.$$
Then the assignment ${\sf U}^{\gamma}_{\sf af}$ defines a $\mathbb{C}$-linear functor 
$${\sf U}^{\gamma}_{\sf af}\colon\mathscr{C}_{[\lambda_{\sf af}]}(V_{\sf af})
\to\mathscr{C}_{[(\lambda+\gamma)_{\sf af}]}(V_{\sf af}),$$
which is a categorical isomorphism with the inverse 
${\sf U}^{-\gamma}_{\sf af}$.
\item
We have an isomorphism
$$\Omega^{-}(\lambda+\gamma)\simeq{\sf U}^{\gamma}_{\sf af}\circ
\Omega^{-}(\lambda)$$
of functors from $\mathscr{C}_{[\lambda]}(V_{\sf sc})$
to $\mathscr{C}_{[(\lambda+\gamma)_{\sf af}]}(V_{\sf af})$.
\end{enumerate}
\end{thm}

\begin{proof}
To prove (1) and (2), it suffices to construct 
a natural isomorphism from
$\Omega^{-}(\lambda+\gamma)$
 to
${\sf U}^{\gamma}_{\sf af}\circ\Omega^{-}(\lambda)$.
For each $\alpha\in\Delta^{+}$, we set
$$(H_{\alpha}^{-})^{*}:=\sum_{\beta\in\Delta^{+}}G^{*}_{\alpha,\beta}H_{\beta}^{-}\in\mathcal{H}^{-}.$$
By using the formula
$\alpha^{*}=\sum_{\beta\in\Pi}c_{\alpha,\beta}H_{\beta}$
and tedious but straightforward calculation,
we obtain
\begin{align*}
\iota^{-}\Big(h_{\sf af}[\gamma]\otimes\mathbf{1}^{-}
+\mathbf{1}_{\sf af}\otimes h^{\gamma}\Big)
=\xi^{\gamma}\otimes\mathbf{1}_{V^{-}}
+\mathbf{1}_{\sf sc}\otimes \zeta^{\gamma}_{-1}\mathbf{1}_{V^{-}},
\end{align*}
where
$$h^{\gamma}:=\sum_{\alpha\in\Delta^{+}}
\gamma(J_{\alpha}^{*})(H_{\alpha}^{-})^{*},\ 
\xi^{\gamma}:=\sum_{\alpha\in\Delta^{+}}
\gamma(J_{\alpha}^{*})\xi(\alpha),\ 
\zeta^{\gamma}:=
\sum_{\alpha\in\Delta^{+}}\gamma(J_{\alpha}^{*})f^{-}_{\sf af}(\alpha).$$
Since we have
a $\mathbb{Z}/2\mathbb{Z}$-graded weak $\mathcal{H}^{-}$-module isomorphism
$$\Big(\mathcal{H}^{-}_{\nu^{-}(\lambda)}\Big)^{h^{\gamma}}\xrightarrow{\simeq}
\mathcal{H}^{-}_{\nu^{-}(\lambda+\gamma)};\,
\ket{\nu^{-}(\lambda)}\mapsto\ket{\nu^{-}(\lambda+\gamma)},$$
in the same way as Theorem \ref{spf_equiv1}, 
we can verify that the natural isomorphism
$$\mathcal{M}\otimes V^{-}\xrightarrow{\simeq}
{\sf U}^{\xi^{\gamma}}_{\sf sc}(\mathcal{M})\otimes
\boldsymbol{\Pi}^{\langle\zeta^{\gamma},\zeta^{\gamma}\rangle}
(V^{-})^{\zeta^{\gamma}_{-1}\mathbf{1}_{V^{-}}}
$$
of $\mathbb{Z}/2\mathbb{Z}$-graded weak $V_{\sf sc}\otimes V^{-}$-modules
for an object $\mathcal{M}$ of $\mathscr{C}_{[\lambda]}(V_{\sf sc})$
gives rise to a natural isomorphism
$\Omega^{-}(\lambda+\gamma)(\mathcal{M})
\simeq{\sf U}^{\gamma}_{\sf af}\circ
\Omega^{-}(\lambda)(\mathcal{M})$
of $\mathbb{Z}/2\mathbb{Z}$-graded weak $V_{\sf af}$-modules.
\end{proof}

\section{Regular cases}

Throughout this section, we assume that $k$ is a positive integer.
In this section, we give additional information about
the representation theory of $\mathcal{C}_{k}(\mf{g})$.

\subsection{Bicommutant of $\mathcal{H}^{+}$}

Let $Q_{l}$ be the even integral sublattice of 
the root lattice $Q$ spanned by long roots.
By the normalization of $(?,?)$ on $\mf{h}^{*}$, we have
$$Q_{l}=\bigoplus_{\alpha\in\Pi}\mathbb{Z}\check{\alpha}
\ \ \left(\check{\alpha}:=\frac{2\alpha}{(\alpha,\alpha)}\right).$$ 
We first recall that
C.\,Dong and Q.\,Wang give a proof of the following fact
(see also \cite[\S4]{dong2009w} for the case of $\mf{g}=\mf{sl}_{2}$).

\begin{prp}[{\cite[Proposition 4.11]{dong2011parafermion}}]
The vertex operator algebra extension 
$$\mathcal{E}:={\sf Com}\Big({\sf Com}\big(\mathcal{H},L_{k}(\mf{g})\big),L_{k}(\mf{g})\Big)$$
of the Heisenberg vertex subalgebra $\mathcal{H}$
generated by $\mf{h}$ in $L_{k}(\mf{g})$ is isomorphic to
the lattice vertex operator algebra associated with
$\sqrt{k}Q_{l}$.
\end{prp}

The next proposition is proved in a similar way as \cite[Proposition 4.11]{dong2011parafermion}
and we omit the detail (see \cite[\S8.2]{CREUTZIG2019396} in the case of $\mf{g}=\mf{sl}_{2}$).

\begin{prp}\label{plus_extension}
The vertex operator superalgebra extension 
$$\mathcal{E}^{+}:={\sf Com}\Big({\sf Com}\big(\mathcal{H}^{+},L_{k}(\mf{g})\otimes V^{+}\big),
L_{k}(\mf{g})\otimes V^{+}\Big)$$
of the Heisenberg vertex subalgebra $\mathcal{H}^{+}$ in $L_{k}(\mf{g})\otimes V^{+}$
is purely even and isomorphic to the lattice vertex operator algebra 
associated with $\sqrt{k+h^{\vee}}Q_{l}$.
Moreover, we have
$$\mathcal{E}^{+}={\sf Com}\Big({\sf Com}\big(\mathcal{H}^{+},
L_{k}(\mf{g})\otimes(V^{+})^{\bar{0}}\big),
L_{k}(\mf{g})\otimes(V^{+})^{\bar{0}}\Big).$$
\end{prp}

\begin{proof}
Note that the even part $(V^{+})^{\bar{0}}$ is 
the lattice vertex operator algebra associated with
the even sublattice 
$$\{\alpha\in L^{+}\,|\,\langle\alpha,\alpha\rangle\in2\mathbb{Z}\}
\simeq
\begin{cases}
2\mathbb{Z} &\text{ if }\mf{g}=\mf{sl}_{2},\\
(\text{the root lattice of type }D_{N}) &\text{ otherwise }
\end{cases}
$$
of $L^{+}$. Then we can apply a straightforward modification of 
the proof of \cite[Proposition 4.11]{dong2011parafermion}.
\end{proof}

\subsection{Bicommutant of $\mathcal{H}^{-}$}

By (\ref{branch2}) and Theorem \ref{spf_equiv2},
we have an isomorphism
\begin{equation*}
\bigoplus_{\gamma\in Q_{\sf sc}}
{\sf U}^{\gamma}_{\sf af}\big(L_{k}(\mf{g})\big)
\otimes
\mathcal{H}^{-}_{\nu^{-}(\gamma)}\simeq 
\iota_{-}^{*}\big(\mathcal{C}_{k}(\mf{g})\otimes V^{-}\big)
\end{equation*}
of $\mathbb{Z}/2\mathbb{Z}$-graded weak $L_{k}(\mf{g})\otimes\mathcal{H}^{-}$-modules.
We first study the structure of
${\sf U}^{\gamma}_{\sf af}\big(L_{k}(\mf{g})\big)$
by the following lemma.

\begin{lem}\label{integrable}
Let $\theta\in P^{\vee}$. Then
$L_{k}(\mf{g})^{\theta}:=\big(L_{k}(\mf{g})\big)^{\theta_{-1}\mathbf{1}_{\sf af}}$
is isomorphic to $L_{k}(\mf{g})$ as a weak $L_{k}(\mf{g})$-module
if and only if $\theta\in Q^{\vee}$.
\end{lem}

\begin{proof}
Define
a Lie algebra automorphism ${\sf t}_{\theta}$ of $\widehat{\mf{g}}=\widehat{\mf{g}}_{B}$ by
\begin{equation*}
{\sf t}_{\theta}\colon\widehat{\mf{g}}\to\widehat{\mf{g}};\ 
X_{\alpha,n}\mapsto X_{\alpha,n+\alpha(\theta)},\ 
H_{n}\mapsto H_{n}+B(H,\theta)K\delta_{n,0},\ K\mapsto K
\end{equation*}
for $\alpha\in\Delta$, $n\in\mathbb{Z}$, and $H\in\mf{h}$.
Then, by a direct computation, we obtain
\begin{equation*}
\displaystyle
Y_{V_{\sf af}}\Big(\Delta\big(\theta_{-1}\mathbf{1}_{\sf af},z\big)X_{-1}\mathbf{1}_{\sf af},z\Big)
=\sum_{n\in\mathbb{Z}}{\sf t}_{\theta}(X_{n})z^{-n-1}
\end{equation*}
for any $X\in\mf{g}$.
By the same argument in \cite[(4-1)]{hosono1991lie},
the automorphism ${\sf t}_{\theta}$ is inner if $\theta\in Q^{\vee}$.
On the other hand, since ${\sf t}_{\theta}$ for $\theta\in P^{\vee}\setminus Q^{\vee}$
corresponds to a non-trivial Dynkin diagram automorphism
and the induced $P^{\vee}/Q^{\vee}$-action on the set of special indices
is simply transitive (see e.g.\,\cite[\S1.2]{minoru2001lectures} for details),
we thus conclude that
the highest weight of $L_{k}(\mf{g})^{\theta}$
coincides with that of $L_{k}(\mf{g})$ 
if and only if $\theta\in Q^{\vee}$.
\end{proof}

\begin{rem}
Following \cite[Proposition 4.1]{hosono1991lie},
we call the automorphism $\mf{t}_{\theta}$ the spectral flow automorphism
of $\widehat{\mf{g}}$ associated with $\theta\in P^{\vee}$.
\end{rem}

We use the next lemma to determine the vertex superalgebra structure
of the bicommutant of $\mathcal{H}^{-}$.

\begin{lem}[cf.\,{\cite[Theorem 4.1]{creutzig2018schur}}]\label{characterization}
Let $(\mathcal{V}=\mathcal{V}^{\bar{0}}\oplus\mathcal{V}^{\bar{1}},
Y,\mathbf{1},\omega)$ be a simple
conformal vertex superalgebra
and $\big(L,\langle?,?\rangle\big)$ a non-degenerate integral lattice.
Assume that there exists an injective $\mathbb{Z}$-linear map
$i\colon L\to\mathcal{V}^{\bar{0}}$ such that
\begin{enumerate}
\item $L_{n}i(h)=\delta_{n,0}i(h)$ for any $h\in L$ and $n\geq0$,
\item $i(h)_{n}i(h')=\delta_{n,1}\langle h,h'\rangle\mathbf{1}$ for any $n\geq0$,
\item we have an isomorphism
$$\bigoplus_{h\in L}(\mathcal{H}_{L})^{i(h)}\simeq\mathcal{V}$$
of (not necessarily $\mathbb{Z}/2\mathbb{Z}$-graded) weak $\mathcal{H}_{L}$-modules,
where $\mathcal{H}_{L}$ is the Heisenberg vertex subalgebra generated by
$i(L)$ in $\mathcal{V}$ and $(\mathcal{H}_{L})^{i(h)}$ is the weak $\mathcal{H}_{L}$-module 
defined in Lemma \ref{LiDelta}.
\end{enumerate}
Then $\mathcal{V}$ is isomorphic to the lattice vertex superalgebra associated with 
the lattice $L$ as a vertex superalgebra.
\end{lem}

\begin{proof}
It follows from a straightforward generalization of \cite[Theorem 3.14]{li1995characterization}
(see \cite[Theorem 6.3.1]{xu1998introduction} for details).
\end{proof}

Now we prove the following proposition.

\begin{prp}\label{minus_extension}
The conformal vertex supralgebra extension
$$\mathcal{E}^{-}:={\sf Com}\Big(
{\sf Com}\big(\mathcal{H}^{-},\mathcal{C}_{k}(\mf{g})\otimes V^{-}\big),
\mathcal{C}_{k}(\mf{g})\otimes V^{-}\Big)$$
of the Heisenberg vertex subalgebra
$\mathcal{H}^{-}$ in $\mathcal{C}_{k}(\mf{g})\otimes V^{-}$ is isomorphic to 
the lattice vertex superalgebra associated with 
$\sqrt{-(k+h^{\vee})}Q_{l}\oplus\mathbb{Z}^{N-\ell}.$
\end{prp}

\begin{proof}
In order to apply Lemma \ref{characterization}, we verify the assumptions in Lemma \ref{characterization}.

We first prove that $\mathcal{E}^{-}$ is simple.
By Proposition \ref{strong_FST} and the regularity of $L_{k}(\mf{g})$,
the vertex superalgebra
$\mathcal{C}_{k}(\mf{g})\otimes V^{-}$
decomposes into a direct sum of simple ordinary
${\sf Com}\big(\mathcal{H}^{-},\mathcal{C}_{k}(\mf{g})\otimes V^{-}\big)$-modules.
Since $V^{-}$ is strongly $L^{-}$-graded
in the sense of \cite[Definition 2.23]{huang2007logarithmic_arXiv},
we can define the contragredient $V^{-}$-module 
in the sense of \cite[Definition 2.35]{huang2007logarithmic_arXiv}.
Then, we can apply the same argument in \cite[Lemma 2.1]{arakawa2017orbifolds}
to our case and conclude that $\mathcal{E}^{-}$ is simple.

Next we verify the conditions (1), (2), and (3) in Lemma \ref{characterization}.
Define a $\mathbb{Z}$-linear map 
$i\colon\sqrt{-(k+h^{\vee})}Q_{l}\oplus\mathbb{Z}^{N-\ell}\to(\mathcal{E}^{-})^{\bar{0}}$ by
\begin{align*}
i\big(\sqrt{-(k+h^{\vee})}\gamma\big)
:=\sum_{\alpha\in\Delta^{+}}(\gamma,\alpha)(H_{\alpha}^{-})^{*},\ \ 
i(\varepsilon_{\beta}):=(H_{\beta}^{-})^{*}
\end{align*}
for $\gamma\in Q_{l}$ and $\beta\in\Delta^{+}\setminus\Pi$,
where $\{\varepsilon_{\beta}\,|\,\beta\in\Delta^{+}\setminus\Pi\}$
is an orthonormal $\mathbb{Z}$-basis of $\mathbb{Z}^{N-\ell}$.
Then the condition (1) is clear and
the condition (2)
follows from (\ref{dualG}) and direct computations.
At last, we verify the condition (3).
It suffices to prove that, for $\gamma\in Q_{\sf sc}$,
the $\mathbb{Z}/2\mathbb{Z}$-graded weak $L_{k}(\mf{g})$-module
${\sf U}^{\gamma}_{\sf af}\big(L_{k}(\mf{g})\big)$ is isomorphic to
$L_{k}(\mf{g})$ if and only if $h^{\gamma}$ lies in the image of the map $i$.
This is verified by Lemma \ref{integrable} and some computations.
We thus complete the proof.
\end{proof}

\subsection{Regularity and Unitarity}

In this subsection, we prove that 
$\mathcal{C}_{k}(\mf{g})$ is regular and has a unitary structure.

We first prove the regularity.
By Proposition \ref{plus_extension}, we obtain
the following coset realization
\begin{equation}\label{condition}
\mathcal{C}_{k}(\mf{g})^{\bar{0}}
={\sf Com}\big(\mathcal{H}^{+},
L_{k}(\mf{g})\otimes(V^{+})^{\bar{0}}\big)
={\sf Com}\big(\mathcal{E}^{+},
L_{k}(\mf{g})\otimes(V^{+})^{\bar{0}}\big).
\end{equation}
Then we obtain the following.

\begin{thm}\label{regular}
Both $\mathcal{C}_{k}(\mf{g})^{\bar{0}}$
and $\mathcal{C}_{k}(\mf{g})$ are regular.
\end{thm}

\begin{proof}
We first prove the regularity of $\mathcal{C}_{k}(\mf{g})^{\bar{0}}$.
By \cite[Theorem 4.5]{abe2004rationality},
it suffices to prove that 
$\mathcal{C}_{k}(\mf{g})^{\bar{0}}$ is rational and $C_{2}$-cofinite.
Then, by (\ref{condition}) and Proposition \ref{plus_extension}, the $C_{2}$-cofiniteness follows from 
\cite[Corollary 2]{miyamoto2015c_2} and
the rationality follows from \cite[Theorem 4.12]{creutzig2018schur}.

Next we prove the regularity of $\mathcal{C}_{k}(\mf{g})$.
By a straightforward generalization of \cite[Theorem 4.5]{abe2004rationality},
it suffices to prove that 
$\mathcal{C}_{k}(\mf{g})$ is rational and $C_{2}$-cofinite.
The rationality of $\mathcal{C}_{k}(\mf{g})$ follows 
from \cite[Proposition 3.6]{dong2012some}.
The $C_{2}$-cofiniteness of $\mathcal{C}_{k}(\mf{g})$ follows from 
\cite[Lemma 2.4]{miyamoto2004modular} and the fact that
$\mathcal{C}_{k}(\mf{g})$ is finitely generated
over the regular vertex operator algebra $\mathcal{C}_{k}(\mf{g})^{\bar{0}}$.
\end{proof}

\begin{cor}\label{ordinary_dec} We have
$$\underline{\mathcal{C}_{k}(\mf{g})\text{-\sf gmod}}
=\bigoplus_{[\lambda]\in P_{\sf sc}^{k}/Q_{\sf sc}}\underline{\mathscr{C}_{[\lambda]}\big(
\mathcal{C}_{k}(\mf{g})\big)},$$
where
$P_{\sf sc}^{k}:=\big\{\mu_{\sf sc}\in (\mf{h}_{\sf sc})^{*}\,\big|\,\mu\in P\big\}.$
\end{cor}

\begin{proof}
The right hand side is a full subcategory of the left hand side
by the definition.
Since the left hand side is a semisimple
abelian category with finitely many (up to isomorphism) simple objects,
it suffices to prove that every simple object in the left hand side
lies in the right hand side.
Let $\mathcal{M}$ be a
simple $\mathbb{Z}/2\mathbb{Z}$-graded ordinary
$\mathcal{C}_{k}(\mf{g})$-module.
Since $\mathcal{M}$ decomposes into a direct sum of finite-dimensional
$L_{0}$-eigenspaces,
the abelian Lie algebra $\mf{h}_{\sf sc}$ acts diagonally on $\mathcal{M}$.
Therfore, by Proposition \ref{converse}, Theorem \ref{MainThm}, and \cite[Theorem 3.1.3]{frenkel1992vertex}, 
there exists $\mu\in P$ such that
$\mathcal{M}$
is a simple object of $\mathscr{C}_{[\mu_{\sf sc}]}\big(\mathcal{C}_{k}(\mf{g})\big)$.
\end{proof}

We next prove the unitarity.

\begin{lem}
There exists an anti-linear involution 
$\phi$ of $\mathcal{C}_{k}(\mf{g})$ such that
$\big(\mathcal{C}_{k}(\mf{g}),\phi\big)$ forms
a \emph{unitary vertex operator superalgebra}
in the sense of \cite[\S2.1]{ai2017unitary}.
\end{lem}

\begin{proof}
By \cite[Theorem 4.7]{dong2014unitary} and \cite[Proposition 2.4 and Theorem 2.9]{ai2017unitary},
both $L_{k}(\mf{g})\otimes V^{+}$ 
and $\mathcal{E}^{+}$
have a unitary structure.
Then, by \cite[Corollary 2.8]{dong2014unitary},
the unitary structure of $L_{k}(\mf{g})\otimes V^{+}$
inherits to $\mathcal{C}_{k}(\mf{g})$.
\end{proof}

We close this subsection with the following result.

\begin{prp}
Every simple $\mathbb{Z}/2\mathbb{Z}$-graded ordinary 
$\mathcal{C}_{k}(\mf{g})$-module appears
as a direct summand of the tensor product
of some simple $\mathbb{Z}/2\mathbb{Z}$-graded ordinary $L_{k}(\mf{g})$-module
and the lattice vertex operator superalgebra $V^{+}$.
In particular, every simple $\mathbb{Z}/2\mathbb{Z}$-graded ordinary
$\mathcal{C}_{k}(\mf{g})$-module is unitarizable.
\end{prp}

\begin{proof}
Let $\mathcal{M}$ be a simple $\mathbb{Z}/2\mathbb{Z}$-graded ordinary 
$\mathcal{C}_{k}(\mf{g})$-module.
By Corollary \ref{ordinary_dec}, we may assume that
$\mathcal{M}$
is a simple object of $\mathscr{C}_{[\mu_{\sf sc}]}\big(\mathcal{C}_{k}(\mf{g})\big)$
for some $\mu\in P$.
Then, by Theorem \ref{MainThm}, we obtain
$\mathcal{M}\simeq\Omega^{+}(\mu)\circ
\Omega^{-}(\mu_{\sf sc})(\mathcal{M})$ as a direct summand of 
$\Omega^{-}(\mu_{\sf sc})(\mathcal{M})\otimes V^{+}$.
Therefore the unitarizability of $\mathcal{M}$ follows from
\cite[Theorem 4.8]{dong2014unitary} and the same discussion as
the proof of \cite[Corollary 2.8]{dong2014unitary}.
\end{proof}

\subsection{Braided monoidal structure}

In this subsection, owing to the recent developments 
in the general theory of vertex superalgebra extensions by
\cite{huang2015braided}, \cite{creutzig2019simple}, and \cite{creutzig2017tensor},
we prove that the category of
$\mathbb{Z}/2\mathbb{Z}$-graded ordinary $\mathcal{C}_{k}(\mf{g})$-modules
carries a braided monoidal category structure.

As a consequence of the theory of vertex tensor categories developed by
Y.-Z.\,Huang and J.\,Lepowsky 
(see \cite{huang1995theory1},
\cite{huang1995theory2},
\cite{huang1995theory3},
\cite{huang1995theory4},
\cite{huang2005differential}, and references therein), we obtain the following.

\begin{prp}\label{evenpart}
The semisimple $\mathbb{C}$-linear abelian category of
purely even $\mathbb{Z}/2\mathbb{Z}$-graded ordinary
$\mathcal{C}_{k}(\mf{g})^{\bar{0}}$-modules
has a braided monoidal category structure
induced by \cite[Theorem 3.7]{huang2005differential}
(see also \cite[Theorem 4.2 and 4.4]{huang1994tensor}).
\end{prp}

\begin{proof}
Since $\mathcal{C}_{k}(\mf{g})^{\bar{0}}$ is regular by Theorem \ref{regular},
all the conditions in \cite[Theorem 3.5]{huang2005differential}
are satisfied.
Therefore it follows from \cite[Theorem 3.7]{huang2005differential}.
\end{proof}

We now arrive at the last result in this paper.

\begin{thm}\label{BTC}
The semisimple $\mathbb{C}$-linear abelian category
$$\mathscr{C}:=\underline{\mathcal{C}_{k}(\mf{g})\text{-\sf gmod}}$$
of $\mathbb{Z}/2\mathbb{Z}$-graded ordinary $\mathcal{C}_{k}(\mf{g})$-modules
has a braided monoidal category structure
induced by \cite[Theorem 3.65]{creutzig2017tensor}.
\end{thm}

\begin{proof}
By \cite[Theorem 3.14]{creutzig2019simple},
the category $\mathscr{C}$
is equivalent to the underlying category 
of the supercategory ${\sf Rep}^{0}\big(\mathcal{C}_{k}(\mf{g})\big)$ defined in
\cite[Definition 2.16]{creutzig2019simple}
(see also \cite[Definition 1.8]{kirillov2002q} for details).
Since these categories carry a natural braided monoidal category structure
by \cite[Theorem 3.65]{creutzig2017tensor} 
(see also \cite[Remark 3.56 and 3.67]{creutzig2017tensor}), we complete the proof.
\end{proof}

\begin{rem}
When $\mf{g}=\mf{sl}_{2}$, Theorem \ref{BTC} is proved by Y.-Z.\,Huang and
A.\,Milas in \cite[Theorem 4.8]{huang2002intertwining}.
Note that Proposition \ref{strong_FST} for $\mf{g}=\mf{sl}_{2}$
(see also \cite[Theorem 5.1]{Ad99}) plays a crucial role 
in their proof of the convergence
and extension property for products of intertwining operators.
See \cite[\S3]{huang2002intertwining} for details.
\end{rem}

\begin{ex}
Assume that $\mf{g}$ is simply-laced and $k=1$.
By Proposition \ref{plus_extension} and 
the isomoprhism $L_{1}(\mf{g})\simeq V_{Q}$ known as the Frenkel--Kac construction, we have
$\mathcal{C}_{1}(\mf{g})\simeq{\sf Com}\big(V_{\sqrt{1+h^{\vee}}Q},V_{Q\oplus\mathbb{Z}^{N}}\big).$
Then this coset vertex operator superalgebra turns out to be isomorphic to the lattice vertex superalgebra associated with
$$Q^{\vee}_{\sf sc}:=\bigoplus_{\alpha\in\Delta^{+}}\mathbb{Z}J_{\alpha},
\ \ \langle J_{\alpha},J_{\beta}\rangle:=(\alpha,\beta)+\delta_{\alpha,\beta}\ 
(\alpha,\beta\in\Delta^{+}).$$
By a direct computation, we can verify that its dual quotient group is given by
$$(Q^{\vee}_{\sf sc})^{*}\!\big/Q^{\vee}_{\sf sc}\simeq
\big(\mathbb{Z}/(1+h^{\vee})\mathbb{Z}\big)^{\ell}.$$
\end{ex}

\appendix

\section{Vertex superalgebras and modules}

\subsection{Universal affine VOAs}\label{Vaff}

Throughout this subsection, $\mf{g}$ stands for 
a finite-dimensional reductive Lie algebra over $\mathbb{C}$.
Let $B\colon \mf{g}\times\mf{g}\to\mathbb{C}$ 
be a symmetric invariant bilinear form and
 $\widehat{\mf{g}}_{B}=\mf{g}\otimes\mathbb{C}[t,t^{-1}]
\oplus\mathbb{C}K$
the corresponding affinization of $\mf{g}$,
{\it i.e.}
the commutation relations are given by
\begin{equation*}
[X_{n},Y_{m}]=[X,Y]_{n+m}+B(X,Y)\,n\delta_{n+m,0}K,\ 
[\widehat{\mf{g}},K]=\{0\}
\end{equation*}
for $X,Y\in\mf{g}$ and $n,m\in\mathbb{Z}$,
where $X_{n}$ stands for $X\otimes t^{n}\in\widehat{\mf{g}}_{B}$.
The next lemma is well-known.

\begin{lem}\label{UnivAff}
Let $\mathbb{C}_{B}$ be a $1$-dimensional 
representation of $\widehat{\mf{g}}_{\geq0}:=\mf{g}\otimes\mathbb{C}[t]\oplus\mathbb{C}K$
defined by $\mf{g}\otimes\mathbb{C}[t].1=\{0\}$ and $K.1=1$.
Then there exists a unique vertex algebra structure on the
induced module
$$V(\mf{g},B):=U(\widehat{\mf{g}}_{B})\otimes_{U(\widehat{\mf{g}}_{\geq0})}
\mathbb{C}_{B}$$
which is strongly generated by
$$X(z):=\sum_{n\in\mathbb{Z}}X_{n}z^{-n-1}$$
for $X\in\mf{g}$, which is
called the \emph{universal affine vertex algebra associated with $(\mf{g},B)$}.
\end{lem}

When $\mf{g}$ is an $n$-dimensional abelian Lie algebra
and $B$ is non-degenerate,
we call $V(\mf{g},B)$ the \emph{Heisenberg vertex algebra of rank $n$}.
It is well-known that Heisenberg vertex algebras of fixed rank
are simple and isomorphic to each other.

\subsection{Coset vertex superlagebras}\label{cvsa}

Let $\mathcal{V}$ be a vertex superalgebra and $\mathcal{W}$ its vertex subsuperalgebra.
Then the $\mathbb{Z}/2\mathbb{Z}$-graded subspace
$${\sf Com}(\mathcal{W},\mathcal{V}):=\{A\in \mathcal{V}\mid [A_{m},B_{n}]=0
\text{ for any }B\in \mathcal{W}\text{ and }m,n\in\mathbb{Z}\}$$
gives a vertex subsuperalgebra of $\mathcal{V}$,
called the coset vertex superalgebra with respect to the pair $(\mathcal{W},\mathcal{V})$.
By the general theory of vertex superalgebras (see e.g.\,\cite{kac97V}), the even linear map
$${\sf Com}(\mathcal{W},\mathcal{V})\otimes \mathcal{W}\to\mathcal{V};\ A\otimes B\to A_{-1}B$$
gives a vertex superalgebra homomorphism.

\subsection{Weak modules}\label{weakmod}

Let $(\mathcal{V}=\mathcal{V}^{\bar{0}}\oplus\mathcal{V}^{\bar{1}},Y,\mathbf{1})$
be a vertex superalgebra.
Let $M=M^{\bar{0}}\oplus M^{\bar{1}}$ be a $\mathbb{Z}/2\mathbb{Z}$-graded vector space and
$$Y_{M}=Y_{M}(?;z)\colon V\rightarrow\End_{\mathbb{C}}(M)[\![z,z^{-1}]\!]$$
a not necessarily $\mathbb{Z}/2\mathbb{Z}$-homogeneous (resp.\,even) linear map.
We write
$$Y_{M}(A;z)=\sum_{n\in\mathbb{Z}}A^{M}_{n}z^{-n-1}\in \End_{\mathbb{C}}(M)[\![z,z^{-1}]\!]$$
for $A\in \mathcal{V}$.
Then the above pair $(M,Y_{M})$ is called a (resp.\,\emph{$\mathbb{Z}/2\mathbb{Z}$-graded})
\emph{weak $\mathcal{V}$-module} if the following conditions hold for any $A\in \mathcal{V}^{\bar{a}}$, $B\in \mathcal{V}^{\bar{b}}$, and $v\in M$:
\begin{description}
\item[Field axiom]
$A^{M}_{n}v=0$ holds for $n\gg0$,
\item[Vacuum axiom]
${\bf 1}^{M}_{n}v=\delta_{n,-1}v$ holds for any $n\in\mathbb{Z}$,
\item[The Borcherds identity]
\begin{equation*}
\begin{array}{ll}
&\displaystyle\sum_{\ell=0}^{\infty}\binom{p}{\ell}(A_{r+\ell}B)^{M}_{p+q-\ell}v\\
&\displaystyle
\hspace{4mm}
=\sum_{\ell=0}^{\infty}(-1)^{\ell}\binom{r}{\ell}
\left(A^{M}_{p+r-\ell}B^{M}_{q+\ell}-(-1)^{r+ab}B^{M}_{q+r-\ell}A^{M}_{p+\ell}\right)v
\end{array}
\end{equation*}
holds for any $p,q,r\in\mathbb{Z}$.
\end{description}

Let $(M^{i},Y_{M^{i}})$ be (resp.\,$\mathbb{Z}/2\mathbb{Z}$-graded) weak $\mathcal{V}$-modules for $i\in\{1,2\}$.
A not necessarily $\mathbb{Z}/2\mathbb{Z}$-homogeneous
(resp.\,even) linear map $f\colon M^{1}\rightarrow M^{2}$ is 
a \emph{morphism of} (resp.\,\emph{$\mathbb{Z}/2\mathbb{Z}$-graded}) \emph{weak $V$-modules} if
$$f\circ A^{M^{1}}_{n}=A^{M^{2}}_{n}\circ f\in\Hom_{\mathbb{C}}(M^{1},M^{2})$$
for any $A\in \mathcal{V}$ and $n\in\mathbb{Z}$.

Let $\mathcal{V}\text{-\sf Mod}$ be 
the $\mathbb{C}$-linear abelian category of weak $\mathcal{V}$-modules
and $\mathcal{V}\text{-\sf gMod}$ the full subcategory\footnote{
We note that the $\mathbb{C}$-linear additive category 
$\mathcal{V}\text{-\sf gMod}$ is not abelian in general.}
of $\mathcal{V}\text{-\sf Mod}$ whose objects
are $\mathbb{Z}/2\mathbb{Z}$-graded
weak $\mathcal{V}$-modules.

\begin{lem}\label{Supercategory}
The categories $\mathcal{V}\text{-\sf Mod}$ and 
$\mathcal{V}\text{-\sf gMod}$ are \emph{$\mathbb{C}$-linear supercategories},
that is, each space of morphisms
is a $\mathbb{Z}/2\mathbb{Z}$-graded vector space over $\mathbb{C}$
and composition of morphisms
induces an even $\mathbb{C}$-linear map.
In addition, the \emph{underlying category}
$\underline{\mathcal{V}\text{-\sf gMod}}$
of the supercategory $\mathcal{V}\text{-\sf gMod}$ in the sense of \cite[Definition 1.1]{brundan2017monoidal} coincides with the $\mathbb{C}$-linear abelian category of
$\mathbb{Z}/2\mathbb{Z}$-graded
weak $\mathcal{V}$-modules.
\end{lem}

\subsection{Ordinary modules}\label{ordinary}

In this subsection, we assume that
$(\mathcal{V}=\mathcal{V}^{\bar{0}}\oplus\mathcal{V}^{\bar{1}},Y,\mathbf{1})$
is a vertex operator superalgebra with respect to a conformal vector 
$\omega\in\mathcal{V}^{\bar{0}}$.
Then a (resp.\,$\mathbb{Z}/2\mathbb{Z}$-graded) weak $\mathcal{V}$-module 
$(M,Y_{M})$ 
is called a (resp.\,\emph{$\mathbb{Z}/2\mathbb{Z}$-graded}) \emph{ordinary
$\mathcal{V}$-module} if $M$ decomposes into
a direct sum
$$M=\bigoplus_{h\in\mathbb{C}}M_{h}$$
of (resp.\,$\mathbb{Z}/2\mathbb{Z}$-homogeneous) finite-dimensional subspaces
$M_{h}:=\{v\in M\,|\,L_{0}v:=\omega^{M}_{1}v=hv\}$
such that $M_{h+n}=\{0\}$ for 
any $h\in\mathbb{C}$ and $n\ll0$.

Similarly to the previous subsection,
we write $\mathcal{V}$-{\sf mod}
(resp.\,$\mathcal{V}\text{-\sf gmod}$)
for the full subsupercategory of 
$\mathcal{V}$-{\sf Mod} whose objects
are (resp.\,$\mathbb{Z}/2\mathbb{Z}$-graded) ordinary $\mathcal{V}$-modules.
Then the underlying category
$\underline{\mathcal{V}\text{-\sf gmod}}$ of 
the $\mathbb{C}$-linear supercategory
$\mathcal{V}\text{-\sf gmod}$ is the $\mathbb{C}$-linear abelian category of
$\mathbb{Z}/2\mathbb{Z}$-graded ordinary $\mathcal{V}$-modules.

\section{Twisted sector}

In this section, we give a remark on the twisted sector of $V_{\sf sc}$
for $k\in\mathbb{C}\setminus\{0,-h^{\vee}\}$.

We set
$$J^{*}:=\frac{1}{2}\sum_{\alpha\in\Delta^{+}}J_{\alpha}^{*}\in\frac{1}{2}P^{\vee}_{\sf sc}.$$
Then the parity involution $\sigma\colon V_{\sf sc}\to V_{\sf sc}$,
which is defined by $\sigma|_{V_{\sf sc}^{\bar{i}}}=(-1)^{i}\id_{V_{\sf sc}^{\bar{i}}}$ for $i\in\{0,1\}$,
satisfies that $\sigma={\sf exp}(2\pi\sqrt{-1}J^{*}_{0})$.

\begin{lem}
The category
of (resp.\,$\mathbb{Z}/2\mathbb{Z}$-graded)
$\mathbb{Z}$-twisted positive energy $V_{\sf sc}$-modules\footnote{
Every $\mathbb{Z}$-twisted positive energy $V_{\sf sc}$-module
is a weak $V_{\sf sc}$-module by the definition.
See \cite[Example 2.14]{de2006finite} for details.} is isomorphic to 
that of (resp.\,$\mathbb{Z}/2\mathbb{Z}$-graded)
$\frac{1}{2}\mathbb{Z}$-twisted positive energy $V_{\sf sc}$-modules.
Here we refer the reader to \cite[Definition 2.21]{de2006finite} for the definition
of $\Gamma$-twisted positive energy modules, where $\Gamma=\mathbb{Z}$ or 
$\frac{1}{2}\mathbb{Z}$.
\end{lem}

\begin{proof}
For a (resp.\,$\mathbb{Z}/2\mathbb{Z}$-graded)
$\mathbb{Z}$-twisted positive energy $V_{\sf sc}$-module $(\mathcal{M},Y_{\mathcal{M}})$,
 we define 
$$Y_{\mathcal{M}}^{\sigma}=Y_{\mathcal{M}}^{\sigma}(?,z)\colon
V_{\sf sc}\to
{\sf End}(\mathcal{M})[\![z^{\frac{1}{2}},z^{-\frac{1}{2}}]\!]$$
by
$Y_{\mathcal{M}}^{\sigma}(v,z):=Y_{\mathcal{M}}\big(\Delta(-J^{*},z)v,z\big)$
for $v\in V_{\sf sc}$.
Then $\mathcal{M}^{\sigma}:=(\mathcal{M},Y_{\mathcal{M}}^{\sigma})$
is a (resp.\,$\mathbb{Z}/2\mathbb{Z}$-graded)
$\frac{1}{2}\mathbb{Z}$-twisted positive energy $V_{\sf sc}$-module and
the assignment $\mathcal{M}\mapsto
\mathcal{M}^{\sigma}$
gives rise to a desired categorical isomorphism.
\end{proof}

\begin{rem}
By \cite[Theorem 2.5]{hosono1991lie},
every $\mathbb{Z}$-twisted (resp.\,$\frac{1}{2}\mathbb{Z}$-twisted) weak $V_{\sf sc}$-module
admits an action of the $\mathcal{N}=2$ superconformal algebra
in the Neveu--Schwarz sector (resp.\,in the Ramond sector)
of central charge $c_{\sf sc}$.
\end{rem}

\bibliographystyle{alpha}

\bibliography{ref}

\end{document}